\documentclass{CSML}

\pdfoutput=1

\usepackage{lastpage}

\def\dOi{13(3:22)2017}
\lmcsheading%
{\dOi}
{1--\pageref{LastPage}}
{}
{}
{Mar.~14, 2016}
{Sep.~13, 2017}
{}

\usepackage{hyperref}\hypersetup{hidelinks}

\usepackage{amssymb, latexsym}
\usepackage{xypic}
\usepackage{tikz}
\usetikzlibrary{arrows.meta}

\DeclareMathOperator{\Pow}{\mathsf{Pow}}
\DeclareMathOperator{\Fin}{\mathsf{Fin}}
\DeclareMathOperator{\PFin}{\mathsf{Fin}^{+}}
\DeclareMathOperator{\cov}{\,\lhd\,}
\DeclareMathOperator{\bcov}{\mathrel{\blacktriangleleft}}
\DeclareMathOperator{\sat}{\mathcal{A}}
\DeclareMathOperator{\amp}{\,\&\,}
\DeclareMathOperator{\imp}{\,\rightarrow\,}
\DeclareMathOperator{\meets}{\between}
\DeclareMathOperator{\subsets}{\subseteq}
\DeclareMathOperator{\Nat}{\mathbb{N}}
\DeclareMathOperator{\Rat}{\mathbb{Q}}
\DeclareMathOperator{\Real}{\mathbb{R}}
\DeclareMathOperator{\UR}{\mathbb{R}^{\mathit{u}}}
\DeclareMathOperator{\NNR}{\mathbb{R}^{\geq0}}
\DeclareMathOperator{\FNNR}{\mathcal{R}^{\geq0}}
\DeclareMathOperator{\FR}{\mathcal{R}}
\DeclareMathOperator{\FUR}{\mathcal{R}^{\mathit{u}}}
\DeclareMathOperator{\PRat}{\mathbb{Q}^{>0}}

\DeclareMathOperator{\URat}{\mathbb{I}}
\DeclareMathOperator{\UInt}{\mathcal{I}[0,1]}
\DeclareMathOperator{\defeqiv}{\stackrel{\textup{def}}{\iff}}
\DeclareMathOperator{\defeql}{\stackrel{\textup{def}}{\  =\  }}

\DeclareMathOperator{\ball}{\mathsf{b}} \DeclareMathOperator{\Ball}{\mathrm{B}}
\DeclareMathOperator{\creg}{\prec\!\!\!\prec\!\!\!\prec}

\DeclareMathOperator{\uc}{\textsf{\textup{uc}}}
\DeclareMathOperator{\mv}{\mathsf{mv}}
\DeclareMathOperator{\rc}{\mathsf{rc}}
\DeclareMathOperator{\wb}{\mathsf{wb}}
\DeclareMathOperator{\Pos}{\mathsf{Pos}}
\DeclareMathOperator{\Sat}{\mathrm{Sat}}
\DeclareMathOperator{\Mod}{\mathcal{M}od}
\DeclareMathOperator{\Pt}{\mathrm{Pt}}

\newcommand{\norm}[1]{\lVert{#1}\rVert}

\DeclareMathOperator{\id}{\mathrm{id}}
\DeclareMathOperator*{\medwedge}{\textstyle{\bigwedge}}

\newcommand{\LKUSpa}{\textbf{\textrm{\textup{LKUSpa}}}}

\newcommand{\GUS}{\textbf{\textrm{\textup{GUS}}}}
\newcommand{\GMS}{\textbf{\textrm{\textup{GMS}}}}
\newcommand{\Ent}{{\mathrm{Rad}}}
\newcommand{\OLKReg}{\textbf{\textrm{\textup{OLKCReg}}}}
\newcommand{\FTop}{\textbf{\textrm{\textup{FTop}}}}

\newcommand{\One}{\textbf{\textrm{\textup{1}}}}
\newcommand{\CZF}{\textmd{\textrm{\textup{CZF}}}}
\newcommand{\REA}{\textmd{\textrm{\textup{REA}}}}
\newcommand{\ACw}{{\ensuremath{\textrm{AC}_{\omega}}}}

\newcommand{\op}{\mathrm{op}}

\theoremstyle{definition}
\newtheorem{notation}[thm]{Notation}
\newtheorem{example}[thm]{Example}

\begin{document}

\title[Localic completion of uniform spaces]{Localic completion of uniform spaces\rsuper*}

\author[T.~Kawai]{Tatsuji Kawai}
\address{
School of Information Science\\
Japan Advanced Institute of Science and Technology\\\newline
1-1 Asahidai\\Nomi\\Ishikawa 923-1292\\Japan\footnote{Author's current
  address: Dipartimento di Matematica, Universit\`{a} di Padova,
             via Trieste 63, 35121 Padova (Italy)}}
\email{tatsuji.kawai@jaist.ac.jp}

\keywords{constructive mathematics; formal topologies; uniform spaces; localic completion; locally compact spaces}
\subjclass{F.4.1}
\titlecomment{{\lsuper*}Some of the results in this paper first appeared in the author's
Ph.D. thesis
\cite{KawaiThesis}}

\begin{abstract}
We extend the notion of localic completion of generalised metric
spaces by Steven Vickers to the setting of generalised uniform spaces.
A generalised uniform space (gus) is a set $X$ equipped with a family
of generalised metrics on $X$, where a generalised metric on $X$ is a
map from $X \times X$ to the upper reals satisfying zero self-distance
law and triangle inequality.

For a symmetric generalised uniform space, the localic completion
lifts its generalised uniform structure to a point-free generalised
uniform structure. This point-free structure induces a complete
generalised uniform structure on the set of formal points of the
localic completion that gives the standard completion of the original gus
with Cauchy filters.

We extend the localic completion to a full and faithful
functor from the category of locally compact uniform spaces into that
of overt locally compact completely regular formal topologies.
Moreover, we give an elementary characterisation of the
cover of the localic completion of a locally compact uniform space
that simplifies the existing characterisation for metric spaces.
These results generalise the corresponding results for metric spaces
by Erik Palmgren.
 
Furthermore, we show that the localic completion
of a symmetric gus is equivalent to the point-free completion of the
uniform formal topology associated with the gus.

We work in Aczel's constructive set theory $\CZF$ with the Regular
Extension Axiom. Some of our results also require Countable Choice.
\end{abstract}
\maketitle
%
\section{Introduction}\label{sec:Intro}
Formal topology \cite{Sambin:intuitionistic_formal_space} is a
predicative presentation of locales, and has been successful in
constructivising many results of the classical topology. However, the
relation between formal topology, which is a point-free notion, and
other constructive point-set approaches to topology is not as simple
as in the classical case.

Classically, the adjunction between the category of topological spaces
and that of locales provides a fundamental tool which relates 
locale theory and point-set topology. Although this adjunction
has been shown to be constructively valid by Aczel
\cite[Theorem 21]{AspectofTopinCZF20063}, it seems to be of little
practical use in Bishop constructive mathematics \cite{Bishop-67}
since we cannot obtain useful correspondence between
point-set notions and point-free notions through this adjunction. For
example, we cannot show that the formal reals and formal Cantor space
are equivalent to their point-set counterparts without employing Fan
theorem \cite[Theorem 3.4]{FormalSpace} \cite[Proposition
4.3]{SchusterGambinoSpatiality}, which is not valid in the recursive
realizability interpretation (see \cite[Chapter 4, Section
7]{ConstMathI}). Since Bishop emphasises the computational aspect
of his mathematics \cite[Appendix B]{Bishop-67}, it seems that Fan
theorem is not acceptable in Bishop constructive
mathematics, the development of which we respect and follow.
Moreover, from the topos theoretic point of view
\cite{SheavesGeometryLogic}, Fan theorem is not desirable as it is not
valid in some Grothendieck toposes \cite[Section 4]{FourmanHyland}.

Instead, Palmgren used Vickers's
notion of localic completion of generalised metric spaces
\cite{LocalicCompletionGenMet} to construct another embedding from the
category of locally compact metric spaces into that of formal
topologies that has more desirable properties, e.g.\ preservation of compactness and local
compactness, and the order of real valued continuous functions \cite{PalmgrenLocalicCompletion}.
This allows us to transfer results in Bishop's theory of metric spaces
to formal topologies.

The main aim of this paper is to further strengthen the connection
between Bishop constructive mathematics and formal topology
by generalising the results by Vickers and Palmgren to the setting of
uniform spaces. The notion of uniform space that we deal with in this
paper is a set equipped with a family of generalised metrics on it, where
a generalised metric is the notion obtained by dropping the symmetry
from that of pseudometric and allowing it to take values in the
non-finite upper reals (Definition \ref{def:GenMet}). 
This structure will be called a generalised uniform space (gus).
When specialised to a family of pseudometrics (i.e.\ finite Dedekind
symmetric generalised metrics),
the notion corresponds to that of uniform space treated in the book by
Bishop \cite[Chapter 4, Problems 17]{Bishop-67}.\footnote{Other
well-known notions of uniform spaces are a set equipped with an
entourage uniformity and a set equipped with a covering uniformity.} 
Hence, the notion of generalised uniform space provides a natural setting to talk
about generalisations of the results by Vickers and Palmgren at the
same time. Our first aim is to extend the notion of localic completion
of generalised metric spaces by Vickers to gus's
(Section \ref{sec:LocComp}). Our second aim is to 
extend Palmgren's results about Bishop metric spaces
\cite{PalmgrenLocalicCompletion} to Bishop uniform spaces by
applying the localic completion of gus's (Section
\ref{sec:FunctEmb}). Our third aim is to relate the localic completion
to the point-free completion of uniform formal topologies by Fox
\cite[Chapter 6]{Fox05} (Section \ref{sec:PFComp}). 

We now summarise our main results.
In Section \ref{sec:LocComp}, we define the localic completion of a
generalised uniform space, and examine its functorial properties and uniform structure.
In particular, we show that the localic completion preserves inhabited
countable products of gus's (Theorem \ref{prop:OmegaProd}).
For a symmetric gus $X$ with a family of generalised metrics $M$, we
show that the localic completion of $X$ is embedded into the
point-free product of the localic completion of each generalised
metric in $M$ as an overt weakly closed subtopology (Theorem \ref{thm:ClosedEmbedding}). This is a
point-free analogue of the usual point-set completion of $X$.
Next, we show that the formal points of the
localic completion of $X$ gives the standard completion of $X$ in
terms of Cauchy filters (Theorem \ref{thm:Completion}).
In Section \ref{sec:FunctEmb}, we specialise the localic completion to
the class of finite Dedekind symmetric gus's, and extend Palmgren's
functorial embedding of locally compact metric spaces to uniform spaces (Theorem \ref{thm:LocalicEmbedding}).
Here, one of our important contributions is a new characterisation of the cover
of the localic completion of a locally compact uniform space
(Proposition \ref{lem:CoverLCM}) that simplifies the previous
characterisation for locally compact metric
spaces by Palmgren \cite[Theorem 4.17]{PalmgrenLocalicCompletion}. 
In Section \ref{sec:PFComp}, we show that the localic completion of a
symmetric gus $X$ is equivalent to the point-free completion of the
uniform formal topology associated with $X$ (Theorem
\ref{thm:EquivLCompPFUComp}).

We work informally in Aczel's constructive set theory $\CZF$ extended
with the Regular Extension Axiom ($\REA$) \cite[Section
5]{Aczel-Rathjen-Note}.  $\REA$ is needed to define the notion of
inductively generated formal topology (see Section
\ref{subsec:InducGenFTop}).
Some of the results in Section \ref{sec:PtLComp} also require Countable
Choice ($\ACw$) (see Remark \ref{rem:Predicativity}).
We assume Palmgren's work \cite{PalmgrenLocalicCompletion} which
is carried out in Bishop's na\"ive set theory \cite{Bishop-67} with some
generalised inductive definitions. It is our understanding
that his work can be carried out in $\CZF + \REA + \ACw$, but we
take extra care not to let our results depend on $\ACw$ implicitly.

\begin{notation}
We fix some notations. For any set $S$, $\Pow(S)$ denotes the class
of subsets of $S$.  Note that since $\CZF$ is predicative, $\Pow(S)$
cannot be shown to be a set unless  $S = \emptyset$. 
$\Fin(S)$ denotes the \emph{set} of finitely enumerable subsets of
$S$, where a set $A$ is finitely enumerable if there exists a
surjection $f \colon \left\{0,\dots,n-1 \right\} \to A$ for some $n
\in \Nat$. $\PFin(S)$ denotes the set of inhabited finitely enumerable subsets of $S$. 
For subsets $U,V \subseteq S$, we define
\begin{align*}
  U \meets V \defeqiv
  \left( \exists a \in S \right) a \in U \cap V.
\end{align*}
Given a relation $r \subseteq S \times T$ between sets $S$ and $T$ and
their subsets $U \subseteq S$ and $D \subseteq T$, we define
\begin{align*}
  r U &\defeql \left\{ b \in T \mid \left( \exists a \in U\right) a
  \mathrel{r} b \right\}, \\
  r^{-} D &\defeql \left\{ a \in S \mid \left( \exists b \in D\right)
    a \mathrel{r} b \right\},\\
  r^{-*} U &\defeql \left\{ b \in T \mid r^{-}\left\{ b \right\}
\subseteq U  \right\}.
\end{align*}
We often write $r^{-}b$ for $r^{-}\left\{ b
\right\}$. The set theoretic complement of a subset $U \subseteq S$
is denoted by $\neg U$, i.e.\
$\neg U
\defeql \left\{ a \in S \mid \neg (a \in U) \right\}$.
\end{notation}

\section{Formal topologies}\label{sec:Prelim}
This section provides background on formal topologies to our main
results in subsequent sections. Our main reference of formal
topologies is \cite{Fox05}, where overt formal topologies are called
open formal topologies. The exposition of the point-free real
numbers in Section \ref{sec:Reals} follows that of Vickers \cite{LocalicCompletionGenMet}.
Nothing in section is essentially new; hence a knowledgeable reader is
advised to skip this section. 

It should be noted that many results on formal topologies are
originally obtained in Martin-L\"of's type theory. So
far, not all of them are available in $\CZF
+ \REA \;( + \ACw)$ (see \cite[Introduction]{vandenBergNID}); however, 
our results do not depend on those results which are only available
in type theory. In any case, we will provide references to the
corresponding results in $\CZF$ in case the original results were
obtained in type theory.

\begin{defi}
  A \emph{formal topology} $\mathcal{S}$ is a triple $(S, \cov, \leq)$
  where $(S, \leq)$ is a preordered set and $\cov$ is a relation
  between $S$ and $\Pow(S)$ such that 
  \[
  \sat U \defeql \left\{ a \in S \mid a \cov U \right\}
  \]
  is a set for each $U \subseteq S$ and that
  \begin{enumerate}
    \item $U \cov U$,
    \item $a \cov U \amp U \cov V \implies a \cov V$,
    \item $a \cov U \amp a \cov V \implies a \cov  U \downarrow V $,
    \item $a \leq b \implies a \cov b$ 
  \end{enumerate}
  for all $a,b \in S$ and $U,V \subseteq S$, where
  \begin{align*}
     U \cov V &\defeqiv \left( \forall a \in U \right) a \cov V,\\
     U \downarrow V &\defeql \left\{ c \in S \mid \left( \exists a \in U
     \right) \left( \exists b \in V \right) c \leq a \amp c \leq
     b\right\}.
  \end{align*}
  We write $a \downarrow U$ for $\left\{ a \right\} \downarrow
  U$ and $U \cov a$ for $U \cov \left\{ a \right\}$.
  The set $S$ is called the \emph{base} of $\mathcal{S}$,
  and the relation $\cov$ is called a \emph{cover} on $S$.
  
  The class $\Sat(\mathcal{S}) = \left\{ \sat U \mid U \in
  \Pow(S) \right\}$ forms a frame with the top $S$, the
  meet $\sat U \wedge \sat V = \sat(U \downarrow V)$ and the join
  $\bigvee_{i \in I}\sat U_i = \sat \bigcup_{i\in I} U_i$ for each $U,V
  \subsets S$ and for each set-indexed family $(U_i)_{i \in I}$ of
  subsets of $S$.

  Formal topologies will be denoted by $\mathcal{S},
  \mathcal{S}',\dots$, and their underlying bases and covers will be denoted
  by $S, S', \dots$, and $\cov, \cov', \dots$ respectively.
\end{defi}

\begin{defi} \label{def:FTM}
  Let $\mathcal{S}$ and $\mathcal{S}'$ be formal topologies. 
  A relation $r \subseteq S \times S'$ is called a \emph{ formal
  topology map} from $\mathcal{S}$ to $\mathcal{S}'$ if 
  \begin{enumerate}[label=(FTM\arabic*)]
    \item\label{FTM1} $S \cov r^{-}S'$,
    \item\label{FTM2} $r^{-}\{a\} \downarrow r^{-}\{b\} \cov r^{-}(a \downarrow'
      b)$,
    \item\label{FTM3} $a \cov' U \implies r^{-}a \cov r^{-}U$
  \end{enumerate}
  for all $a, b \in S'$ and $U \subseteq S'$. 
  The class $\mathrm{Hom}(\mathcal{S}, \mathcal{S}')$ of formal topology maps
  from $\mathcal{S}$ to $\mathcal{S}'$ is ordered by
  \[
    r \leq s \defeqiv \left( \forall a \in S' \right) r^{-}a \cov
    s^{-}a.
  \]
  Two formal topology maps $r,s : \mathcal{S} \to \mathcal{S}'$
  are defined to be \emph{equal}, denoted by $r = s$, if $r \leq s$
  and $s \leq r$.

  A formal topology map $r \colon \mathcal{S} \to \mathcal{S}'$
  bijectively corresponds to a frame homomorphism $F_{r} \colon
  \Sat(\mathcal{S}') \to \Sat(\mathcal{S})$ between the associated
  frames in such a way that  $F_{r}(\sat' U) = \sat r^{-} U$.
\end{defi}

The formal topologies and formal topology maps between them form a
category $\FTop$.  The composition of two formal topology maps is the
composition of the underlying relations of these maps.  The identity
morphism on a formal topology is the identity relation on its base.

The formal topology $\One \defeql \left(\left\{ * \right\}, \in, =
\right)$ is a terminal object in $\FTop$. A formal point of a 
formal topology $\mathcal{S}$ is a formal topology map $r \colon \One \to
\mathcal{S}$. An equivalent definition is the following.
\begin{defi}\label{def:Pt}
  Let $\mathcal{S}$ be a formal topology. A subset $\alpha \subseteq S$
  is a \emph{formal point} of $\mathcal{S}$ if
  \begin{enumerate}[label=({P}\arabic*)]
    \item\label{P1} $ S \meets \alpha$,
    \item\label{P2} $ a, b \in \alpha \implies \alpha \meets (a \downarrow
      b)$,
    \item\label{P3} $ a \in \alpha \amp a \cov U \implies \alpha \meets U$
  \end{enumerate}
  for all $a, b \in S$ and $U \subseteq S$.
  The class of formal points of $\mathcal{S}$ is denoted by
  $\Pt(\mathcal{S})$. Predicatively, we cannot assume that
  $\Pt(\mathcal{S})$ is a set.\footnote{For example, for any set $S$,
    the class of formal points of a formal topology $(\Fin(S),
    \cov, \leq)$, where 
    $A \leq B \defeqiv B \subseteq A$ and 
    $A \cov \mathcal{U} \defeqiv \left( \exists B \in
    \mathcal{U} \right) A \leq B$,
    is in one-one correspondence with $\Pow(S)$.}
\end{defi}
Impredicatively,  the formal points $\Pt(\mathcal{S})$
form a topological space with the topology generated
by the basic opens of the form
\[
  a_{*} \defeql \left\{ \alpha \in \Pt(\mathcal{S}) \mid a \in
\alpha \right\}
\]
for each $a \in S$. If $r \colon \mathcal{S} \to \mathcal{S}'$ is a
formal topology map, then the function $\Pt(r) \colon
\Pt(\mathcal{S}) \to \Pt(\mathcal{S}')$ given by
\begin{equation*}
  \Pt(r)(\alpha)  \defeql r \alpha
\end{equation*}
for each $\alpha \in \Pt(\mathcal{S})$ is a well defined continuous
function with respect to the topologies on $\Pt(\mathcal{S})$ and
$\Pt(\mathcal{S}')$. Impredicatively, the operation $\Pt(-)$ is the
right adjoint of the adjunction between the category of
topological spaces and formal topologies established by Aczel \cite{AspectofTopinCZF20063}.

From Section \ref{sec:LocComp} onward, we mainly work with overt
formal topologies, i.e.\ formal topologies equipped with a positivity predicate.
\begin{defi}
  Let $\mathcal{S}$ be a formal topology.
  A subset $V \subseteq S$ is said to be \emph{splitting} if 
  \[
  a \in V \amp  a \cov U \implies V \meets U
  \]
  for all $a \in S$ and $U \subseteq S$.
  A \emph{positivity predicate} (or just a positivity) on a formal
  topology $\mathcal{S}$ is a splitting subset $\Pos \subseteq S$
  that satisfies
  \begin{equation}\label{Pos}
    \tag{Pos} a \cov \left\{ x \in S \mid x = a \amp \Pos(a) \right\}
  \end{equation}
  for all $a \in S$, where we write $\Pos(a)$ for $a \in \Pos$.
  A formal topology is \emph{overt} if it is equipped with a
  positivity predicate.
\end{defi}
Let $\mathcal{S}$ be a formal topology.  By condition \eqref{Pos},
a positivity predicate on $\mathcal{S}$, if it exists, is the largest
splitting subset of $\mathcal{S}$. Thus, a formal topology admits at most
one positivity predicate.

\subsection{Inductively generated formal
topologies}\label{subsec:InducGenFTop} The notion of inductively
generated formal topology by Coquand et al.\ \cite{Coquand200371}
allows us to define a formal topology by a small set of axioms.

\begin{defi}
Let $S$ be a set. An \emph{axiom-set} on $S$ is a pair $(I,C)$, where
$(I(a))_{a \in S}$ is a family of sets indexed by  $S$, and $C$ is a
family $(C(a,i))_{a \in S, i \in I(a)}$ of subsets  of $S$
indexed by $\sum_{a \in S}I(a)$. 
\end{defi}

The following result, which was obtained in Martin-L\"of's type theory, is also
valid in $\CZF + \REA$ \cite[Section 6]{AspectofTopinCZF20063}.
\begin{thm}[{\cite[Theorem 3.3]{Coquand200371}}] \label{thm:IIFTopi}
  Let $(S ,\leq)$ be a preordered set, and let $(I,C)$ be an axiom-set on
  $S$. Then, there exists a cover $\cov_{I,C}$ inductively generated
  by the following rules:
  \begin{gather*}
    \frac{a \in U}{a \cov_{I,C} U} \text{ (reflexivity)};
    \quad
    \frac{a \leq b \ \  b \cov_{I,C} U}{a \cov_{I,C} U} 
    \text{ ($\leq$-left)}; \\
    \frac{a \leq b \ \  i \in I(b) \ \ \left\{ a \right\} \downarrow
    C(b,i) \cov_{I,C} U}{a \cov_{I,C} U} \text{ ($\leq$-infinity)}.
  \end{gather*}
  The relation $\cov_{I,C}$ is the least cover on $S$ which satisfies
  ($\leq$-left) and $a \cov_{I,C} C(a,i)$ for each $a \in S$ and $i
  \in I(a)$, and  $\cov_{I,C}$ is called the cover inductively
  generated by $(I,C)$. 
\end{thm}
A formal topology $\mathcal{S}=(S, \cov, \leq)$ is
\emph{inductively generated} if it is equipped with an axiom-set $(I,C)$
on $S$ such that $\cov = \cov_{I,C}$.

Localised axiom-sets are particularly convenient to work with.
\begin{defi}
  Let $(S,\leq)$ be a preordered set, and 
  let $(I,C)$ be an axiom-set on  $S$.  We say that $(I,C)$ is
  \emph{localised} with respect to $(S, \leq)$ if
\begin{align*}
  a \leq b \implies
  \left( \forall i \in I(b) \right)\left( \exists j \in I(a) \right)
  C(a,j) \subseteq a \downarrow  C(b,i)
\end{align*}
for all $a,b \in S$.
We often leave implicit the preorder with respect to which the axiom-set is
localised. In that case, the context will always make it clear what is left implicit.
\end{defi}
For an axiom-set $(I,C)$ which is localised with respect to a preorder $(S,
\leq)$, we can replace ($\leq$-infinity) rule in Theorem
\ref{thm:IIFTopi} with (infinity):
\[
    \frac{i \in I(a) \ \ 
    C(a,i) \cov_{I,C} U}{a \cov_{I,C} U} \text{ (infinity)}.
\]

Thus, if $\mathcal{S} = (S,\cov,\leq)$ is the formal topology
inductively generated by a localised axiom-set $(I,C)$,
 then for each $U \subseteq S$, the set
  \begin{align*}
    \sat U \defeql \left\{ a \in S \mid a \cov U \right\}
  \end{align*}
  is the least subset of $S$ such that
  \begin{enumerate}
    \item $U \subseteq \sat U$,
    \item $a \leq b \amp b \in \sat U \implies a \in \sat U$,
    \item $C(a,i) \subseteq \sat U \implies a \in \sat U$
  \end{enumerate}
for all $a,b \in S$ and $i \in I(a)$. 

Therefore, we have the following induction principle:
let $(I,C)$ be a localised axiom-set with respect to a preorder $(S,\leq)$. Then,
for any subset $U \subseteq S$ and a predicate $\Phi$  on $S$,
if
\begin{enumerate}[label=({ID}\arabic*)]
  \item\label{ID1} $\displaystyle{\frac{a \in U}{\Phi(a)}}$,
  \item\label{ID2} $\displaystyle{\frac{a \leq b \quad \Phi(b)}{\Phi(a)}}$,
  \item\label{ID3} $\displaystyle{\frac{i \in I(a)
    \quad \left( \forall c \in C(a,i) \right)
    \Phi(c)}{\Phi(a)}}$
\end{enumerate}
for all $a,b \in S$, then $a \cov_{I,C} U \implies
\Phi(a)$ for all $a \in S$. 
An application of the above principle is called a \emph{proof by
induction on the cover $\cov_{I,C}$}.

\begin{rem}\label{rem:IndGenFTop}
In Definition \ref{def:FTM} of formal topology map,
if the formal topology $\mathcal{S}'$ is inductively generated by an
axiom-set $(I,C)$ on $S'$, then the condition \ref{FTM3} is equivalent to
the following conditions under the condition \ref{FTM2}.
\begin{enumerate}[label=\!\!({FTM3}\alph*)]
    \item\label{FTM3a} $a \leq' b \implies r^{-} a \cov r^{-} b$,
    \item\label{FTM3b} $r^{-} a \cov r^{-}C(a,i)$
\end{enumerate}
for all $a,b \in S'$ and $i \in I(a)$.

Similarly, in Definition \ref{def:Pt} of formal point,
if the formal topology $\mathcal{S}$ is inductively generated by an
axiom-set $(I,C)$ on $S$, then the condition \ref{P3} is equivalent to
the following conditions:
\begin{enumerate}[label=\!\!({P3}\alph*)]
  \item $a \leq b \amp a \in \alpha \implies b \in \alpha$,
  \item $a \in \alpha  \implies \alpha \meets C(a,i)$
\end{enumerate}
for all $a, b \in S$ and $i \in I(a)$.
See Fox \cite[Section 4.1.2]{Fox05} for further details.

In the same setting as above, a subset $V \subseteq S$ is splitting if
and only if the following two conditions hold:
\begin{enumerate}[label=\upshape({Spl}\arabic*)]
  \item\label{Spl1} $a \in V \amp a \leq b \implies b \in V$,
  \item\label{Spl2} $a \leq b \amp a \in V \amp i \in I(b)
        \implies V \meets (a \downarrow C(b,i))$.
\end{enumerate}
Furthermore, if the axiom-set $(I,C)$ is localised, the condition \ref{Spl2} can be
replaced by 
\begin{enumerate}[label=\upshape({Spl}\arabic*'),start=2]
  \item\label{Spl2'} $a \in V \amp i \in I(a)
        \implies V \meets  C(a,i)$.
\end{enumerate}
\end{rem}

\subsection{Products and pullbacks} \label{sec:Limit}
Inductively generated formal topologies are closed under arbitrary
limits. We recall the constructions of products and pullbacks.
\subsubsection{Products}\label{sec:FTopiProduct}
Following Vickers
\cite{Vickers05someconstructive}, we define a product of a
set-indexed family of inductively generated formal topologies as follows.
Let $(\mathcal{S}_i)_{i \in I}$ be a family of inductively generated
formal topologies, each member of which is of the form $ \mathcal{S}_i= (S_i, \cov_i,
\leq_i)$, and let $(K_i, C_i)$ be the axiom-set which generates
$\mathcal{S}_i$. Define a preorder $(S_{\Pi}, \leq_{\Pi})$ by
\begin{align*}
  S_{\Pi} &\defeql \Fin\left(\sum_{i \in I}S_i\right), \\
  A \leq_{\Pi} B &\defeqiv \left( \forall (i,b) \in B\right) \left(
  \exists (j,a) \in A \right) i = j \amp a \leq_i b
\end{align*}
for all $A,B \in S_{\Pi}$.
The axiom-set on $(S_{\Pi},\leq_{\Pi})$ is given by
\begin{enumerate}[label=({S}\arabic*)]
  \item\label{S1} $S_{\Pi} \cov_{\Pi} \left\{ \left\{ (i, a) \right\} \in
    S_{\Pi} \mid a \in S_i\right\}$ for each $i \in I$,
  \item\label{S2} $ \left\{ (i,a), (i, b) \right\} \cov_{\Pi} \left\{
    \left\{ (i,c) \right\} \in S_{\Pi} \mid c \leq_i a \amp c \leq_i b
    \right\}$ for each $i \in I$ and $a,b \in S_i$,
  \item\label{S3} $\left\{ (i,a) \right\} \cov_{\Pi} \left\{ \left\{ (i, b)
    \right\} \in S_{\Pi} \mid b \in C_i(a, k) \right\}$  for each
    $i \in I$, $a \in S_i$, and $k \in K_i(a)$.
\end{enumerate}

Let $\prod_{i \in I} \mathcal{S}_i = (S_{\Pi}, \cov_{\Pi},
\leq_{\Pi})$ be
the formal topology inductively generated by the above axiom-set.
For each $i \in I$, the projection $p_i \colon \prod_{i \in I} \mathcal{S}_{i} \to  \mathcal{S}_i$ is defined by
\[
A  \mathrel{p_i}   a \defeqiv A = \left\{ (i, a) \right\}
\]
for all $A \in S_{\Pi}$ and $a \in S_i$.

Given any family $(r_i \colon \mathcal{S} \to \mathcal{S}_i)_{i \in I}$
of formal topology maps, we have a unique formal topology map $r
\colon \mathcal{S} \to \prod_{i \in I} \mathcal{S}_{i}$ such that
$r_{i} = p_{i} \circ r$ for each $i \in I$. The map $r$ is defined by
\[
a \mathrel{r}  A \defeqiv \left( \forall (i, b) \in A \right) a \cov
r_{i}^{-}  b
\]
for each $a \in S$ and $A \in S_{\Pi}$.

\subsubsection{Binary products}\label{sec:FTopiBiProduct}
A binary product of a pair of inductively generated formal topologies
admits a simple construction. Given two inductively generated formal
topologies $\mathcal{S} = (S, \cov_S, \leq_S)$ and $\mathcal{T} =
(T,\cov_T,\leq_T)$ generated by axiom-sets $(I,C)$ and $(J,D)$
respectively, their product $\mathcal{S} \times \mathcal{T}$ is an
inductively generated formal topology with the preorder $(S \times T,
\leq)$  defined by
\[
(a,b) \leq (a',b') \defeqiv a \leq_S a' \amp b \leq_{T} b'
\]
and the axiom-set $(K,E)$ on $(S \times T, \leq)$  defined by
\begin{align*}
   K\left( (a,b) \right) &\defeql I(a) + J(b),\\
   E\left( (a,b), (0, i)\right) &\defeql C(a,i) \times \left\{ b \right\},\\
   E\left( (a,b), (1, j)\right) &\defeql \left\{ a \right\} \times D(b,j).
\end{align*}
The projection $p_{\mathcal{S}} \colon \mathcal{S} \times \mathcal{T} \to
\mathcal{S}$ is given by
\begin{align*}
  (a,b) \mathrel{p_{\mathrel{S}}} a'
  \defeqiv (a,b) \cov_{K,E} \left\{ a' \right\} \times T
\end{align*}
for each $a,a'\in S$ and $b \in T$, and the other projection is
similarly defined.

Given two formal topology maps $r \colon \mathcal{S}' \to
\mathcal{S}$ and  $s \colon \mathcal{S}' \to \mathcal{T}$, the canonical
map $\langle r,s \rangle \colon \mathcal{S}' \to \mathcal{S} \times
\mathcal{T}$ is given by
\begin{align*}
  c \mathrel{\langle r,s \rangle} \left(a,b\right) \defeqiv c \cov'
  r^{-} a \amp c \cov' s^{-} b
\end{align*}
for each $c \in S'$, $a \in S$, and $b \in T$.

\subsubsection{Pullbacks}\label{sec:FTopiPullback}
Given inductively
generated formal topologies $\mathcal{S}_1$ and $\mathcal{S}_2$
generated by axiom-sets $(I_1,C_1)$ and $(I_2,C_2)$ respectively, a
pullback $\mathcal{S}_1 \times_{\mathcal{T}} \mathcal{S}_2$
of formal topology maps $r \colon \mathcal{S}_1 \to \mathcal{T}$, $s
\colon
\mathcal{S}_2 \to \mathcal{T}$ is generated by the axiom-set of
the product $\mathcal{S}_1 \times \mathcal{S}_2$ together with the following
additional axioms:
\begin{align*}
  (a,b) &\cov S_1 \times s^{-}c \qquad(a \mathrel{r} c \amp c \in T),\\
  (a,b) &\cov r^{-}c \times S_2 \qquad(b \mathrel{s} c \amp c \in T).
\end{align*}
The restrictions of the projections
$p_i \colon \mathcal{S}_1 \times \mathcal{S}_2  \to \mathcal{S}_i \;
(i = 1,2)$ to $\mathcal{S}_1 \times_{\mathcal{T}} \mathcal{S}_2$ form
a pullback square.

\subsection{Subtopologies}\label{sec:SubSpa}
\begin{defi}
  A \emph{subtopology} of a formal topology $\mathcal{S}$ is a formal
  topology $\mathcal{S}' = (S, \cov', \leq)$ where $\cov'$ is a
  cover on $S$ and $(S, \leq)$ is the underlying preorder of $\mathcal{S}$
  such that  $\sat U \subseteq \sat' U$ for each $U \subseteq S$. If
  $\mathcal{S}'$ is a subtopology of $\mathcal{S}$, we write
  $\mathcal{S}'\sqsubseteq \mathcal{S}$.

  Given a formal topology map $r \colon \mathcal{S} \to \mathcal{S}'$, the
  relation $\cov_r \subseteq S' \times \Pow(S')$ given by 
  \[
    a \cov_r U \defeqiv a \in r^{-*}\sat r^{-}U
  \]
  is a cover on $S'$. The formal topology $\mathcal{S}_r  = (S', \cov_r, \leq')$ is
  called the \emph{image} of $\mathcal{S}$ under $r$. 
  A formal topology map $r \colon \mathcal{S} \to \mathcal{S}'$ is 
  an \emph{embedding} if $r$ restricts to an isomorphism between
  $\mathcal{S}$ and its image $\mathcal{S}_{r}$, and $r$ is called a
  \emph{surjection} if its image is $\mathcal{S}'$.
  It can be shown that $r \colon \mathcal{S} \to \mathcal{S}'$
  is an embedding if and only if
  \[
  a \cov r^{-}r^{-*}\sat\left\{ a \right\}
  \]
  for all $a \in S$. See Fox \cite[Proposition 3.5.2]{Fox05}.
\end{defi}
  By the condition \ref{FTM3} for a formal topology map, we have
  $\mathcal{S}_{r} \sqsubseteq \mathcal{S}'$ for any formal topology
  map $r \colon \mathcal{S} \to \mathcal{S}'$.  If $\mathcal{S}'$ is a
  subtopology of $\mathcal{S} = (S,\cov,\leq)$, then the identity
  relation $\id_{S}$ on $S$ is an embedding $\id_{S}
  \colon \mathcal{S}' \to \mathcal{S}$. Hence the notion of embedding
  is essentially equivalent to that of subtopology.

  The following is well known. We omit an easy proof.
\begin{lem}\label{lem:ImgOvert}
  Let $\mathcal{S}$ be an overt formal topology with a positivity
  $\Pos$, and let $r \colon \mathcal{S} \to \mathcal{S}'$ be a formal
  topology map. Then, the image $\mathcal{S}_r$ of $\mathcal{S}$
  under $r$ is overt with the positivity 
  \[
  r \Pos = \left\{ a \in S' \mid
  \left( \exists b \in \Pos \right) b \mathrel{r} a \right\}.
  \] 
\end{lem}

The notion of  weakly closed subtopology due to Bunge and Funk
\cite{BungeFunkLowerPowerLoc} is particularly relevant to us (see
Vickers \cite{SublocFTop} and Fox \cite{Fox05} for treatments in
formal topology). For the reason of predicativity, we only consider
overt weakly closed subtopologies of inductively generated formal
topologies.
\begin{defi}\label{def:WeaklyClosed}
  Let $\mathcal{S}$ be an inductively generated formal topology, and
  let $V \subseteq S$ be a splitting subset of $\mathcal{S}$.
  The overt \emph{weakly closed subtopology} of $\mathcal{S}$
  determined by $V$, denoted by $\mathcal{S}_{V}$, is inductively generated by the axioms
  of $\mathcal{S}$ together with the following extra axioms:
  \[
    a \cov_{V} V \cap \left\{ a \right\}
  \]
  for each $a \in S$. In this case, $\mathcal{S}_{V}$ is overt with
  positivity $V$.
\end{defi}
Clearly, $\mathcal{S}_{V}$ is the
largest subtopology of $\mathcal{S}$ with positivity $V$. Moreover,
we have $\Pt(\mathcal{S}_{V}) = \left\{ \alpha \in
\Pt(\mathcal{S}) \mid \alpha \subseteq V \right\}$ by Remark
\ref{rem:IndGenFTop}, and $\Pt(\mathcal{S}_{V})$ forms a closed
subspace of $\Pt(\mathcal{S})$.

\subsection{Upper reals and Dedekind reals}\label{sec:UpperReals}
We introduce the notions of upper reals and (possibly non-finite) Dedekind reals,
which we use for the values of distance functions of metric spaces.
These real numbers are defined as models of propositional
geometric theories.

\subsubsection{Propositional geometric theories}
Along with the notion of inductively generated formal topology,
propositional geometric theories (i.e.\ presentations of frames by
generators and relations \cite[Chapter 2]{vickers1989topology})
provide us with equivalent, but more logical descriptions of formal
topologies.  Vickers \cite{Vickers06CompactnessLocFTop} gives a good
exposition of the connection between frame presentations and
inductively generated formal topologies (see also \cite[Chapter
4]{Fox05}).

\begin{defi}
Let $G$ be a set, whose elements are called \emph{propositional
symbols} or \emph{generators}. A \emph{propositional geometric theory}
$T$ over $G$
is a set $R_{T} \subseteq \Fin(G) \times \Pow(\Fin(G))$ of
\emph{axioms}. An axiom $(P, \left\{ P_{i} \mid i \in I \right\}) \in
R_{T}$ is usually denoted by
\[
\medwedge P \vdash \bigvee_{i \in I} \medwedge P_i
\]
or
\[
  p_{0} \wedge \cdots \wedge p_{n-1} \vdash \bigvee_{i \in I}
  p^{i}_{0} \wedge \cdots \wedge p^{i}_{n_{i}-1}.
\]
We write $\top$ for $\medwedge \emptyset$.  
Henceforth, propositional geometric theories will be simply called
geometric theories. A geometric theory over
propositional symbols $G$ with axioms $R$ will be denoted by a pair
$(G,R)$.
\end{defi}

Every geometric theory $T = (G,R)$ determines a formal topology
$\mathcal{S}_{T} = (S_{T}, \cov_{T}, \leq_{T})$ in the following way.
The base $S_{T}$ is $\Fin(G)$ ordered by $A \leq_{T} B
\defeqiv B \subseteq A$. The cover $\cov_{T}$ is
generated by the following axioms:
\[
  P \cov_{T} \left\{ P_{i} \mid i \in I \right\}
\]
for each $\medwedge P \vdash \bigvee_{i \in I} \medwedge
P_{i} \in R$.
The topology $\mathcal{S}_{T}$ represents a frame \emph{presented}
by the theory in the following sense:
\begin{enumerate}
  \item there exists a function $\iota_{T} \colon G \to
    \Sat(\mathcal{S_{T}})$, defined by $\iota_{T}(p) = \sat_{T} \left\{
    p \right\}$, which preserves all the axioms in $R$, i.e.\
    \[
     \iota_{T}(p_{0}) \downarrow \cdots \downarrow \iota_{T}(p_{n-1})
     \cov_{T} \bigcup_{i \in I} \iota_{T}(p^{i}_{0}) \downarrow \cdots
     \downarrow \iota_{T}(p^{i}_{n_{i}-1})
    \]
    for each $
     p_{0} \wedge \cdots \wedge p_{n-1} \vdash \bigvee_{i \in I}
     p^{i}_{0} \wedge \cdots \wedge p^{i}_{n_{i}-1} \in R$.
  \item for any frame $X$ and a function $f \colon G \to Y$ which
    preserves all the axioms in $R$, there exists a unique
    frame homomorphism $F \colon \Sat(\mathcal{S}_{T}) \to X$ such
    that $F \circ \iota_{T} = f$, where  $F$ is given by
    \[
      F(\sat U) \defeql \bigvee_{P \in U} \medwedge_{p \in P}
      f(p)
    \]
    for each $U \subseteq S_{T}$.
\end{enumerate}
In particular, each formal topology map $r \colon \mathcal{S} \to
\mathcal{S}_{T}$ is determined by a function $f \colon G \to \Pow(S)$ such that
\[
  f(p_{0}) \downarrow \cdots \downarrow f(p_{n-1}) \cov \bigcup_{i
  \in I} f(p^{i}_{0}) \downarrow \cdots \downarrow  f(p^{i}_{n_{i}-1})
\]
for each
$
  p_{0} \wedge \cdots \wedge p_{n-1} \vdash \bigvee_{i \in I}
  p^{i}_{0} \wedge \cdots \wedge p^{i}_{n_{i}-1} \in R
$.
In this case, we have
\[
  a \cov r^{-}\left\{ P \right\} \iff \left(\forall p \in P  \right) a \cov
  f(p)
\]
for each $P \in \Fin(G)$.

\begin{rem}
  Adding axioms to a geometric theory $T$ amounts to defining a
  subtopology $\mathcal{S}'$ of $\mathcal{S}_{T}$. Indeed, the
  identity function on the propositional symbols $G$ of $T$ gives
  rise to the canonical subspace inclusion from $\mathcal{S}'$ to
  $\mathcal{S}_{T}$ represented by the identity relation on
  $\Fin(G)$.
\end{rem}

\subsubsection{Models of theories}
A \emph{model} of a geometric theory $T = (G,R)$ is a subset $m
\subseteq G$ such that
\[
  P \subseteq m \implies  \left( \exists i \in I \right) P_{i}
  \subseteq m
\]
for each axiom $\medwedge P \vdash \bigvee_{i \in I} \medwedge
P_{i} \in R$. Let $\Mod(T)$ denote the class of models of $T$.
There exists a bijective correspondence between the models of
$T$ and the formal points of $\mathcal{S}_{T}$:
\begin{align*}
  m &\mapsto \Fin(m) \colon \Mod(T) \to \Pt(\mathcal{S}_{T}),\\
  \alpha &\mapsto \left\{ p \in G \mid \left\{ p \right\} \in m
\right\} \colon \Pt(\mathcal{S}_{T}) \to \Mod(T).
\end{align*}
The class $\Mod(T)$ of models gives rise to a topological space,
the topology of which is generated by the subbasics of the form
\[
  a_{*} \defeql \left\{ m \in \Mod(T) \mid a \in m \right\}.
\]
The above bijective correspondence between $\Mod(T)$ and
$\Pt(\mathcal{S}_{T})$ determines homeomorphisms between
the associated spaces of $\Mod(T)$ and $\Pt(\mathcal{S}_{T})$.

\subsubsection{Real numbers}\label{sec:Reals}
Let $\PRat$ denote the set of positive rational numbers.
An \emph{upper real} is a model of the geometric theory $T_{u}$ over
$\PRat$ with the following axioms:
\begin{align*}
  q &\vdash q' &&(q \leq q'),\\
  q &\vdash \bigvee_{q' < q} q'.
\end{align*}
That is, an upper real is a subset $U \subseteq \PRat$ such that
\[
  q \in U \iff \left( \exists q' < q \right) q' \in  U.
\]
The class of the upper reals will be denoted by $\UR$, and the
formal topology determined by the theory of the upper reals
will be denoted by $\FUR$.

Non-negative rational numbers are embedded into $\UR$ by $r \mapsto
\left\{ q \in \PRat \mid r < q\right\}$, which we simply write as
$r$.
The orders on $\UR$ is defined by
\begin{align*}
  U \leq V &\defeqiv V \subseteq U,\\
  U < V &\defeqiv \left( \exists r \in \PRat \right) U + r \leq V.
\end{align*}
Note that for $U \in \UR$ and $q \in \PRat$, we have $U < q \iff q \in U$.

An upper real is \emph{finite} if it is inhabited, i.e.\ if it is a
model of the theory $T_{u}$ with the extra axiom:
\[
  \top \vdash \bigvee_{q \in \PRat} q.
\]

A \emph{non-finite Dedekind real} is a model of the geometric theory
$T_{D}$ over $G_D \defeql \Rat + \Rat$, elements of which will be
denoted by $(p, + \infty) \defeql (0,p)$ and $(- \infty, q) \defeql
(1,q)$. The axioms of $T_{D}$ are the following:
\begin{align*}
  (-\infty,q) &\vdash (-\infty,q') &&(q \leq q') \\
  (-\infty,q) &\vdash \bigvee_{ q' < q} (-\infty,q')  \\
  (p, +\infty) &\vdash (p',+\infty) &&(p' \leq p) \\
  (p, +\infty) &\vdash \bigvee_{p'>p} (p', +\infty)  \\
  (p, +\infty) \wedge (-\infty, q) &\vdash \bigvee
  \left\{ (p, +\infty) \wedge (-\infty, q) \mid p < q \right\} \\
  \top &\vdash (p, +\infty) \vee (-\infty, q) && (p < q) 
\end{align*}
A non-finite Dedekind real $m$ is equivalent to a pair $(L,U)$ of
possibly empty lower cut $L$ and upper cut $U$ with the following
correspondence:
\[
  L \defeql \left\{ p \in \Rat \mid (p, + \infty ) \in m \right\},
  \qquad 
  U \defeql \left\{ q \in \Rat \mid (-\infty, q) \in m \right\}.
\]

The orders and additions on the non-finite Dedekind reals are defined by adding 
the conditions dual to those of the upper reals to the lower cuts.
For example, if $(L,U)$ and $(L',U')$ are non-finite Dedekind reals $(L,U) \leq
(L',U') \defeqiv L \subseteq L' \amp U' \subseteq U$.

A \emph{non-negative non-finite Dedekind real} is a model of
the theory $T_{D}$ extended with the following axioms:
\[
  (-\infty, q) \vdash \bigvee \left\{ (-\infty, q) \mid 0 < q \right\}.
\]
The class of non-negative non-finite Dedekind reals will be denoted
by $\NNR$, and the formal topology determined by the theory of 
non-negative non-finite Dedekind reals will be denoted by $\FNNR$.

A non-finite Dedekind real is \emph{finite} (or just a Dedekind real)
if both its lower and upper cuts are inhabited, i.e.\ if it is a
model of the theory $T_{D}$ with the extra axioms:
\[
  \top \vdash \bigvee_{q \in \Rat} (-\infty, q), \qquad
  \top \vdash \bigvee_{p \in \Rat} (p, +\infty).
\]
Non-negative Dedekind reals are defined similarly and its collection
will be denoted by $\Real^{0+}$.

The embedding $q \in \PRat \mapsto (-\infty, q) \in G_{D}$ of
generators gives rise to a formal topology map $\iota_{D} \colon \FNNR
\to \FUR$. It is not hard to see that $\iota_{D}$
is a monomorphism.\footnote{Hints: use the last three axioms of the theory
  $T_{D}$.} 
The morphism $\iota_{D}$ restricts to a morphism between formal topologies
associated with the theories of finite upper
reals and non-negative Dedekind reals.

If $T_{1} = (G_{1},R_1)$ and $T_{2} = (G_{2},R_2)$ are geometric
theories, then the product $\mathcal{S}_{T_{1}} \times
\mathcal{S}_{T_{2}}$ is presented by generators $G_{1} + G_{2}$ and
axioms
\[
  (k,p_{0}) \wedge \cdots \wedge (k, p_{n-1}) \vdash \bigvee_{i \in I}
  (k, p^{i}_{0}) \wedge \cdots \wedge (k,p^{i}_{n_{i}-1})
\]
for each axiom 
$
  p_{0} \wedge \cdots \wedge p_{n-1} \vdash \bigvee_{i \in I}
  p^{i}_{0} \wedge \cdots \wedge p^{i}_{n_{i}-1}  \in R_{k} \; (k \in \left\{ 0,1 \right\})$.
The injections of generators $p \mapsto (k,p) \;(p \in G_k)$ give
rise to the projections $p_{k} \colon \mathcal{S}_{T_{1}} \times
 \mathcal{S}_{T_2} \to \mathcal{S}_{T_{k}}$.

  The formal addition $ + \colon \FUR \times \FUR \to \FUR$ on $\FUR$
  is determined  by a function $f_{+} \colon \PRat \to \Pow(\Fin(\PRat +
  \PRat))$ given by
  \[
    f_{+}(q) \defeql \left\{  \left\{ (0, q_{1}), (1, q_{2})\right\} \mid
    q_{1} + q_{2} < q \right\},
  \]
  which induces the addition on $\UR$ by
  \[
    U + V \defeql \left\{ q_{1} + q_{2} \mid q_{1} \in U \amp
    q_{2} \in V  \right\}.
  \]
  Similarly, the addition $ + \colon \FNNR \times \FNNR \to \FNNR$ on
  $\FNNR$ is determined by a function $g_{+} \colon G_{D} \to \Pow(\Fin(G_{D} + G_{D}))$
  given by
  \begin{align*}
    g_{+}(-\infty, q) &\defeql \left\{ \left\{ (0, (-\infty, q_{1}),
    (1, (-\infty, q_{2})\right\} \mid
    q_{1} + q_{2} < q \right\}, \\
    g_{+}(p, +\infty) &\defeql \left\{ \left\{ (0, (p_{1}, +\infty),
    (1, (p_{2}, +\infty)\right\} \mid
    p < p_{1} + p_{2} \right\},
  \end{align*}
  which induces the addition on $\NNR$ by
  \[
    (L_{1}, U_{1}) + (L_{2}, U_{2}) \defeql (L_{1} + L_{2}, U_{1} +
    U_{2}).
  \]
  It is easy to see that the following diagram commutes: 
  \[
  \xymatrix@C+=4em{
      \FNNR \times \FNNR \ar[r]^{\iota_D \times \iota_D} \ar[d]_-{+} 
      & \FUR \times \FUR \ar[d]^-{+} \\
     \FNNR \ar[r]^{\iota_D} &  \FUR
      }
  \]

    The formal order $\leq_{u}$ on $\FUR$ is a subtopology of $\FUR
    \times \FUR$ defined by adding the following axioms to the theory
    of
    $\FUR \times \FUR$:
    \[
      (1,q) \vdash (0,q).
    \]
    Similarly, the order $\leq_{D}$ on $\FNNR$ is a subtopology of
    $\FNNR \times \FNNR$ defined by adding the following axioms to the
    theory of
    $\FNNR \times \FNNR$:
    \[
      (0,(p,+\infty)) \wedge (1,(-\infty,q))
      \vdash \bigvee \left\{ (0,(p,+\infty)) \wedge (1,(-\infty,q))
      \mid p < q \right\}.
    \]
  Note that those subtopologies induce the orders on $\UR$ and
  $\NNR$ that have been defined before.
  The morphism $\iota_D \times \iota_D \colon \FNNR \times \FNNR \to
  \FUR \times \FUR$ restricts to a morphism $\iota_D \times \iota_D
  \colon {\leq_{D}} \to {\leq_{u}}$, which we denote
  by the same symbol by an abuse of notation. It is easy to see that the
  following diagram commutes and that it is a pullback diagram:
  \[
  \xymatrix@C+=4em@M+=.5em{
    \leq_{D} \ar[r]^{\iota_D \times \iota_D} \ar@{^{(}->}[d]
      &  \leq_{u} \ar@{^{(}->}[d] \\
     \FNNR \times \FNNR \ar[r]^{\iota_D \times \iota_D} & \FUR \times
     \FUR
      }
  \]

\section{Localic completion of generalised uniform spaces}\label{sec:LocComp}
\subsection{Generalised uniform spaces}
We first recall the notion of generalised metric space by Vickers
\cite{LocalicCompletionGenMet}.
In Bishop's theory of metric space \cite{Bishop-67}, a
distance function takes values in the non-negative Dedekind reals.
Vickers allowed three generalisations to the usual development of
metric spaces: first, the distances need not be finite. Second, the
values are taken in the upper reals (not necessarily in the Dedekind
reals). Third, the distances are not assumed to be symmetric.

\begin{defi}[{\cite[Definition 3.4]{LocalicCompletionGenMet}}]\label{def:GenMet}
  A \emph{generalised metric} on a set $X$ is a function $d \colon X \times X
  \to \UR$ such that 
  \begin{enumerate}
    \item $d(x,x) = 0$,
    \item $d(x,z) \leq d(x,y) + d(y,z)$
  \end{enumerate}
  for all $x,y,z \in X$. A generalised metric $d$ is \emph{finite} if
  $d$ takes values in the finite upper reals and $d$ is called
  \emph{Dedekind} if $d$ factors uniquely through the non-negative
  (non-finite) Dedekind reals. A generalised metric $d$ is
  \emph{symmetric} if
  \[
    d(x,y) = d(y,x)
  \]
  for all $x,y \in X$.

  If $d$ and $\rho$ are generalised metrics on a set $X$, we define
  \begin{equation}\label{eq:OrderGM}
    d \leq \rho \defeqiv \left( \forall x, x' \in X \right) d(x,x')
    \leq \rho(x,x').
  \end{equation}
  
  If $A$ is a finitely indexed set $\left\{ d_{0},\dots,d_{n-1}
\right\}$ of generalised metrics on a set $X$, then the function
  $d_{A} \colon X \times X \to \UR$ defined by
  \[
    d_{A}(x,x') \defeql \sup\left\{ d_{i}(x,x') \mid i < n \right\}
  \]
  is again a generalised metric on $X$.

  A set equipped with a generalised metric
  is called a \emph{generalised metric space} (abbreviated as gms).
  A \emph{homomorphism} from a gms $(X,d)$ to another gms $(Y,\rho)$
  is a function $f \colon X \to Y$ such that 
  \[
    \rho(f(x),f(x')) \leq d(x,x')
  \]
  for each $x,x' \in X$.
Generalised metric spaces and homomorphisms between them form a
category $\GMS$.
\end{defi}

The following seems to be the most natural generalisation of the notion of gms.
\begin{defi}
  A \emph{generalised uniform space} (abbreviated as gus) is a set $X$
  equipped with a set $M$ of generalised metrics on it, where $M$ 
  is inhabited and closed under binary sups with respect to the order
  $\leq$ given by \eqref{eq:OrderGM}.
  A \emph{homomorphism} from a gus $(X,M)$ to a gus
  $(Y,N)$ is a function $f \colon X \to Y$ such that
  \begin{align*}
    \left( \forall \rho \in N \right) \left( \exists d \in M \right)
    \left[ \left( \forall x,x' \in X \right) \rho(f(x),f(x')) \leq
    d(x,x') \right].
  \end{align*}
\end{defi}

A gus $(X,M)$ is called \emph{finite} (Dedekind, symmetric) if
each $d \in M$ is finite (respectively, Dedekind, symmetric).

\begin{rem}\label{rem:BishUnif}
  Bishop \cite[Chapter 4, Problems 17]{Bishop-67} defined a uniform space as a pair $(X,M)$ of
  set $X$ and a set $M$ of pseudometrics on $X$, where $M$ is not assumed 
  be closed under binary sups.\footnote{
    Bishop \cite[Chapter 4, Problems
    17]{Bishop-67} does not even impose 
    inhabitedness on $M$. We decided to include the condition for
    a smooth development of localic completions of
    generalised uniform spaces. By doing so, we can also incorporate
  the theory of generalised metric spaces into that of generalised
uniform spaces more naturally; for example, with our definition of
gus, the inclusion of the category of gms's into that of gus's
preserves finite limits.}
    However, we can equip $X$ with a new set
   $M'$ of pseudometrics on $X$ given by $M' = \left\{ d_{A} \mid A \in
   \PFin(M) \right\}$, which is uniformly isomorphic
    to $(X,M)$ in the sense of \cite[Chapter 4, Problem 17]{Bishop-67}.
  Hence, our assumption on $M$ is compatible with Bishop's approach.
\end{rem}

The generalised uniform spaces and homomorphisms between them form a
category $\GUS$. The category of generalised metric spaces $\GMS$ 
can be embedded into $\GUS$ by $(X,d) \mapsto (X,\left\{d
\right\})$.
Obviously, homomorphisms between gms's become homomorphisms between the
corresponding gus's. We usually identify $(X,\left\{ d \right\})$
with $(X,d)$.

\begin{example}
  Any function $f \colon X \to \Real$ from a set $X$ to the finite
  Dedekind reals $\Real$ determines a pseudometric on $X$ by%
  \[
    d_{f}(x,x') \defeql |f(x) - f(x')|.
  \]
  Hence, any subset of the set $\mathbb{F}(X,\Real)$ of real valued
  functions on $X$ determines a uniform structure on $X$ (see Remark
  \ref{rem:BishUnif}).
  Important examples are
  \begin{itemize}
    \item the set of pointwise continuous functions from a metric
      space to $\Real$;
    \item the set of uniform continuous functions from a compact metric
      space to $\Real$;
    \item the set of continuous functions from a locally compact metric
      space to $\Real$.%
      \footnote{
  For the notions of local compactness and continuous functions 
  used in this examples, see Definition \ref{def:LKUSpa}. }
  \end{itemize}
  Those are examples of \emph{function spaces} \cite[Chapter 3,
  Definition 8]{Bishop-67}, the notion which has gained renewed
  interest in recent years (see
  \cite{BridgesFunctionSpace,IshiharaFunctionSpace,PetrakisUrysohn}).

  Similar examples are obtained when we consider a seminorm
  on a linear space (e.g.\ over $\Real$), where a
  seminorm on a linear space $V$ is a non-negative real valued
  function $\norm{-} \colon V \to \Real^{0+}$ such that
  \begin{align*}
    \norm{av} &= |a|\norm{v}, & \norm{v + w} &\leq \norm{v} +
    \norm{w} && \tag*{($a \in \Real, \; v,w \in V$)}
  \end{align*}
  Any seminorm $\norm{-} \colon V \to \Real^{0+}$ determines a pseudometric
  on $V$ by
  \[
    d(v,w) \defeql \norm{v - w}.
  \]
  Hence, any subset of the collection of seminorms on a linear space
  $V$ determines a uniform structure on $V$. \emph{Locally convex spaces}
  are particularly important examples of the structure of this
  kind, where a locally convex space is a pair $(V,\mathfrak{N})$
  of a linear space $V$ and a set of seminorms $\mathfrak{N}$ on
  $V$ such that for each $\norm{-}_{1},\norm{-}_{2} \in \mathfrak{N}$ and $c \in \PRat$,
  and for each seminorm $\norm{-}$ on $V$, 
  \[
    \left[ \left( \forall v \in V \right) \norm{v} \leq c(\norm{v}_{1} +
    \norm{v}_{2}) \right] \implies
    \norm{-} \in \mathfrak{N}.
  \]
  A simple example of a locally convex space is the ring
  $C(X)$ of real valued continuous functions on a locally compact
  metric space $X$ with the locally convex structure generated by
  the seminorms
  \begin{align*}
    \left\{ \norm{-}_{K} \mid \text{$K \subseteq X$ compact
    subset}\right\},
  \end{align*}
  where $\norm{f}_{K} \defeql \sup \left\{ \lvert f(x) \rvert \mid x \in K \right\}$.
  See \cite[Chapter 9, Section 5]{Bishop-67} for details on locally
  convex spaces.
\end{example}

\begin{example}[{Non-symmetric gus}]
  Lifting the restriction of symmetry allows us to metrise more spaces. 
  If $(X,d)$ is a generalised metric space, then we can define a generalised
  metric $d_{L}$ on $\Fin(X)$ given by
  \[
    d_{L}(A,B) \defeql \sup_{a \in A} \inf_{b \in B}d(a,b),
  \]
  which is called the \emph{lower metric}.  This construction, together
  with the upper and the Hausdorff generalised metric, is treated in
  detail by Vickers \cite{LocCompPowerLocales}.
  These constructions can be naturally extended to 
  generalised uniform spaces. It is interesting to see how much of
  the results in \cite{LocCompPowerLocales} can be carried over to the
  setting of generalised uniformly spaces.
\end{example}

\begin{example}[{Domain theory \cite[Section 5]{LocalicCompletionGenMet}}]
  \label{eg:Domain}
  The notion of generalised metric space (and generalised uniform
  space) and its localic completion to be defined in Section
  \ref{sec:LocComp} provide a common generalisation of 
  the theory of metric space and that of domain (in the sense of domain
  theory \cite{abramsky1994domain}).

  Let $(P,\leq)$ be a poset. Then, we can define a generalised metric on
  $P$ by
  \[
    d(x,y) \defeql \left\{ q \in \PRat \mid x \leq y \vee 1 < q
    \right\}.
  \]
  A more elaborate example is the rationals $\Rat$ with the following
  generalised metric:
  \[
    d(x,y) \defeql \left\{ q \in \PRat \mid x < y + q \right\}.
  \]
  See also Remark \ref{rem:IdealCompletion}.
\end{example}

\subsection{Localic completions}\label{subsec:LocComp}
Given a gus $(X,M)$,  we define the set
\[
  \Ent(M) \defeql M \times \PRat
\]
of generalised radii parameterised by $M$ and a set
\[
  U_X \defeql \Ent(M) \times X
\]
of generalised formal balls.
We write $\ball_d(x,\varepsilon)$ for the element $\left( \left(d,
\varepsilon  \right), x\right) \in U_X$. 

Define an order $\leq_{X}$ and a transitive relation $<_{X}$ on $U_X$ by
\begin{align*}
  \ball_d(y,\delta) \leq_{X} \ball_{\rho}(x,\varepsilon) \defeqiv
  \rho \leq d \amp \rho(x,y) + \delta \leq \varepsilon,\\
  \ball_d(y,\delta) <_{X} \ball_{\rho}(x,\varepsilon) \defeqiv
   \rho \leq d \amp \rho(x,y) + \delta < \varepsilon.
\end{align*}
We extend the relations  $\leq_X$ and $<_X$ to the subsets of
$U_X$ by 
\[
  U \leq_X V \defeqiv \left( \forall a \in U \right) \left( \exists b
  \in V \right) a \leq_X b
\]
for all $U,V \subseteq U_X$, and similarly for $<_X$.

The \emph{localic completion} of a gus $X = (X,M)$ is a formal
topology
\begin{align*}
  \mathcal{U}(X) = (U_{X}, \cov_X, \leq_{X}),
\end{align*}
where $\cov_{X}$ is inductively
generated by the axiom-set on $U_{X}$ consisting of the
following axioms:
\begin{enumerate}[label=({U}\arabic*)]
  \item\label{U1} $a \cov_{X} \left\{ b \in U_X \mid b <_X a \right\}$;
  \item\label{U2} $a \cov_{X} \mathcal{C}_{d}^{\varepsilon}\;$ for each
    $\left(d, \varepsilon \right) \in \Ent(M)$
\end{enumerate}
for each $a \in S$, where we define 
\begin{equation}\label{eq:UniCov}
  \mathcal{C}_{d}^{\varepsilon}  \defeql \left\{ \ball_d(x, \varepsilon)
  \in U_X \mid x \in X \right\}.
\end{equation}

\begin{rem}
For a generalised metric space $(X,d)$, the localic completion of
the gus $(X,\left\{ d \right\})$ is equivalent to the localic completion
$\mathcal{M}(X,d)$ of the gms $(X,d)$ by Vickers
\cite{LocalicCompletionGenMet}, and we use the notation
$\mathcal{M}(X,d)$ to denote  $\mathcal{U}(X,\left\{ d \right\})$.
\end{rem}
\begin{rem}\label{rem:IdealCompletion}
One of the important aspects of the localic completion is that the
notion captures important constructions on symmetric uniform spaces
and domains in a single setting. For example, Vickers showed that the
localic completion of the generalised metric associated with a
poset $(P,\leq)$ given in Example \ref{eg:Domain} represents the
Scott topology on $P$ \cite[Proposition 5.6]{LocalicCompletionGenMet},
and the localic completion of the generalised metric on $\Rat$ in the
same example represents the topology determined by the geometric
theory of the lower cuts (the dual notion of the finite upper
reals defined in Section \ref{sec:Reals}) \cite[Proposition
5.7]{LocalicCompletionGenMet}.

What the localic completion means for the usual symmetric case will be
studied in Section \ref{sec:SymmGUS}. We will not pursue the domain
theoretic aspect of localic completions further in this paper.
\end{rem}

\begin{lem}\label{lem:U2EquivU2'}
  The axioms of the form \ref{U2} are equivalent to the
  following axioms:
\begin{enumerate}[label=({U}\arabic*'),start=2]
  \item\label{U2'} $a \cov_{X} \mathcal{C}_{d}^{\varepsilon}
    \downarrow a$ for each $\left(d, \varepsilon \right) \in \Ent(M)$,
\end{enumerate}
that is, together with \ref{U1}, they generate the same cover on $U_X$.
\end{lem}
\proof
  Obvious.
\qed
Note that the axiom-set consisting of axioms of the forms
\ref{U1} and \ref{U2'} is localised with respect to $\leq_{X}$.

For each $\ball_{d}(x,\varepsilon) \in U_X$, we write
$\ball_d(x,\varepsilon)_{*}$ or $\Ball_d(x,\varepsilon)$  to denote
the \emph{open ball} corresponding to $\ball_{d}(x,\varepsilon)$,
i.e.\
\begin{align*}
  \ball_d(x,\varepsilon)_{*}
  \defeql \Ball_d(x,\varepsilon)
   \defeql \left\{ x' \in X \mid d(x,x') < \varepsilon \right\}.
\end{align*}
We extend the notation $(-)_{*}$ to the subsets of $U_X$ by
  $
  U_{*} \defeql \bigcup_{a \in U} a_{*}.
  $

Dually, each $x \in X$ is associated with the set $\Diamond x$ of
\emph{open neighbourhoods} of $x$, namely 
\[
\Diamond x \defeql \left\{ a \in U_X \mid x \in
a_{*}\right\}.
\]

\begin{lem}\label{lem:misc}
  Let $(X,M)$ be a gus. Then
  \begin{enumerate}
    \item\label{lem:misc1} $a' \leq_X a <_X b \leq_X b' \implies a' <_X b'$,
    \item\label{lem:misc2} $a <_X b \implies \left( \exists c \in U_X \right) a <_X  c <_X
      b$,
    \item\label{lem:misc3} $a \leq_X b \implies a_{*} \subseteq b_{*}$
  \end{enumerate}
  for all $a,a',b,b' \in U_X$.
\end{lem}
\proof
 Straightforward.
\qed

\begin{rem}\label{rem:DifferenceBallWise}
  The converse of Lemma \ref{lem:misc} \eqref{lem:misc3} need
  not hold.
For example, consider the unit interval
$(\left[0,1\right], d)$ of $\Real$, where $d$ denotes the usual metric on
$\left[0,1\right]$. We have $\Ball_{d}(1,3) \subseteq \Ball_{d}(1,2)$, but
$\ball_{d}(1,3) \leq_{[0,1]} \ball_{d}(1,2)$ is false.
\end{rem}

\begin{prop}\label{prop:LCompOvrt}
  For any gus $X$, its localic completion
  $\mathcal{U}(X)$ is overt, and the base $U_X$
  is the positivity of $\mathcal{U}(X)$.
\end{prop}
\proof
  Straightforward.
\qed

\subsection{Functoriality of localic completions}
Given a homomorphism $f \colon (X,M) \to (Y,N)$ of gus's, define a
relation $r_{f} \subseteq U_X \times U_Y$ by
\begin{equation}\label{eq:GUSFTop}
  \ball_{d}(x,\varepsilon) \mathrel{r_{f}} \ball_{\rho}(y,\delta)
  \defeqiv d \mathrel{\omega_{f}} \rho \amp \ball_{\rho}(f(x),
  \varepsilon) <_Y \ball_{\rho}(y,\delta),
\end{equation}
where $d \mathrel{\omega_{f}} \rho \defeqiv \left( \forall x,x' \in X
\right) \rho(f(x),f(x')) \leq d(x,x')$.
\begin{prop}
  The localic completion extends to a functor from $\GUS$ to
  $\FTop$.
\end{prop}
\proof
  We must show that the assignment $f \mapsto r_{f}$ is functorial. 
  First, we show that $r_{f}$ is a formal topology map for any
  homomorphism $f \colon (X,M) \to (Y,N)$. We show that $r_{f}$
  satisfies \ref{FTM2}. The other properties of the morphism are easy to prove.

  \ref{FTM2} Suppose that $\ball_{d}(x,\varepsilon) \in
  r_{f}^{-} \ball_{\rho_{1}}(y_{1}, \delta_{1}) \downarrow  r_{f}^{-}
  \ball_{\rho_{2}}(y_{2}, \delta_{2})$. Then,
  \begin{align*}
    d \mathrel{\omega_{f}} \rho_{1} &\amp
    \ball_{\rho_{1}}(f(x), \varepsilon) <_{Y} \ball_{\rho_1}(y_{1},
    \delta_{1}), \\
    d \mathrel{\omega_{f}} \rho_{2} &\amp
    \ball_{\rho_{2}}(f(x), \varepsilon) <_{Y} \ball_{\rho_2}(y_{2}, \delta_{2}).
  \end{align*}
  Let $\rho = \sup\left\{ \rho_{1}, \rho_{2} \right\}$, and choose
  $\theta \in \PRat$ such that $\ball_{\rho_{1}}(f(x),\varepsilon +
  \theta ) <_{Y} \ball_{\rho_{1}}(y_1, \delta_{1})$  and 
  $\ball_{\rho_{2}}(f(x),\varepsilon +
  \theta ) <_{Y} \ball_{\rho_{2}}(y_2, \delta_{2})$.
  Then, $\ball_{\rho}(f(x),\varepsilon +
  \theta ) \in  \ball_{\rho_{1}}(y_1, \delta_{1}) \downarrow
  \ball_{\rho_{2}}(y_2, \delta_{2})$. Moreover,
  $\rho(f(x),f(x')) \leq d(x,x')$ for all $x,x' \in X$, so that $d
  \mathrel{\omega_{f}} \rho$. Thus, $\ball_{d}(x,\varepsilon) \mathrel{r_{f}}
  \ball_{\rho}(f(x), \varepsilon + \theta)$, from which  \ref{FTM2}
  follows.

  Next we show that the assignment is functorial. 
  First, for any gus $(X,M)$, we have $\ball_{d}(x,\varepsilon)
  \mathrel{r_{\id_{X}}} \ball_{\rho}(y,\delta) \iff \ball_{d}(x,\varepsilon)
  <_X \ball_{\rho}(y,\delta)$ so that $r_{\id_{X}} =
  \id_{\mathcal{U}(X)}$ as formal topology maps.
  Second, let $f \colon (X,M) \to (Y,N)$ and $g \colon (Y,N) \to
  (Z,L)$ be homomorphisms. Let $a = \ball_{d}(x,\varepsilon) \in U_X$
  and $c = \ball_{\rho}(z,\xi) \in U_{Z}$, and
  suppose that $\ball_{d}(x,\varepsilon) \mathrel{r_{g \circ f}}
  \ball_{\rho}(z,\xi)$. Then, $d \mathrel{\omega_{g \circ f}} \rho \amp
  \ball_{\rho}(g(f(x)), \varepsilon) < \ball_{\rho}(z,\xi)$.
  Since $f$ and $g$ are homomorphisms, there exist $\rho' \in N$
  and $d' \in M$ such that $d' \mathrel{\omega_{f}} \rho'
  \mathrel{\omega_{g}} \rho$. Choose $\theta \in \PRat$ such that
  $\ball_{\rho}(g(f(x)), \varepsilon + \theta) <_Z
  \ball_{\rho}(z,\xi)$. Let $\ball_{d^{*}}(x',\varepsilon') \in
  a \downarrow \mathcal{C}^{\theta}_{d'}$. Then,
  we have $\ball_{d^{*}}(x',\varepsilon') \mathrel{r_{f}}
  \ball_{\rho'}(f(x'), \varepsilon' + \theta) \mathrel{r_{g}}
  \ball_{\rho}(z,\xi)$. Thus $r_{g \circ f}^{-} c \cov_{X}
  (r_{g} \circ r_{f})^{-}c$. We also easily have
  $(r_{g} \circ r_{f})^{-}c \cov_{X} r_{g \circ f}^{-} c$. Hence
  $r_{g \circ f} = r_{g} \circ r_{f}$.
\qed
Let us denote this functor by $\overline{({-})} \colon \GUS \to
\FTop$.

The binary product of gus's $(X,M)$ and $(Y,N)$ is defined by
\[
  X \times Y = (X \times Y, M \times N),
\]
where an element $(d,\rho) \in M \times N$ is regarded as a
generalised metric on $X \times Y$ defined by
\[
  (d, \rho)((x,y),(x',y')) \defeql \max\left\{ d(x,x'), \rho(y,y')
  \right\}.
\]
Binary sups in $M \times N$ is defined coordinate-wise. The projections
$\pi_{X} \colon X \times Y \to X$ and $\pi_{Y} \colon X \times Y \to Y$ 
are obviously homomorphisms of gus's. A terminal gus is $(\left\{ *
\right\}, \left\{ d_{*} \right\})$, which is a one-point set
with discrete metric $d_{*}$.

\begin{prop}\label{prop:Prod}
  The functor $\overline{(-)} \colon \GUS \to \FTop$ preserves
  finite products.
\end{prop}
\proof
The proof of the preservation of a terminal object is the same as that
of generalised metric spaces \cite[Proposition 5.3]{LocalicCompletionGenMet}.

We sketch the proof of preservation of binary products of gus's,
which is analogous to the corresponding fact about gms's
\cite[Theorem 5.4]{LocalicCompletionGenMet}.
Given gus's $(X,M)$ and $(Y,N)$, the functor sends the projection
$\pi_{X} \colon X \times Y \to X$ to a formal topology map
$\overline{\pi_{X}} \colon \mathcal{U}(X \times Y) \to
\mathcal{U}(X)$, which can be easily shown to be equal to a
relation $r_{X} \subseteq U_{X \times Y} \times U_{X}$ defined by
\[
  \ball_{(d,\rho)}((x,y), \xi) \mathrel{r_{X}}
  \ball_{d'}(x',\varepsilon') \defeqiv \ball_{d}(x,\xi) <_X
  \ball_{d'}(x',\varepsilon').
\]
There is also a relation $r_{Y} \subseteq U_{X \times Y} \times
U_{Y}$, similarly defined as $r_{X}$,  which is equal to
$\overline{\pi_{Y}}$. 
On the other hand, we define a relation 
$r \subseteq \left( U_{X} \times U_{Y} \right) \times U_{X \times Y}$ by
\[
  \left( \ball_{d}(x,\varepsilon), \ball_{\rho}(y,\delta) \right)
  \mathrel{r} \ball_{(d',\rho')}( (x',y'), \xi)
  \defeqiv 
  \ball_{d}(x,\varepsilon) <_{X} \ball_{d'}( x', \xi)
  \amp
  \ball_{\rho}(y,\delta) <_{Y} \ball_{\rho'}( y', \xi),
\]
which can be easily shown to be a formal topology map.

Then, it is straightforward to show that $\langle r_{X}, r_{Y} \rangle \circ
r = \id_{\mathcal{U}(X) \times \mathcal{U}(Y)}$ and 
$r \circ \langle r_{X}, r_{Y} \rangle  = \id_{\mathcal{U}(X \times
Y)}$. 
\qed

The product of a set-indexed family $(X_{i})_{i \in I}$ of gus's,
each member of which is of the form $X_i = (X_i, M_i)$, consists of the cartesian product
$\prod_{i \in I}X_i$  and the set
\[
  M_{\Pi} \defeql \PFin\left(\sum_{i \in I}M_i\right)
\]
of generalised metrics on $\prod_{i \in I}X_i$, where we identify each member
$A \in M_{\Pi}$ with a generalised metric on $\prod_{i \in I}X_i$
defined by
\[
  A(f,g) \defeql \sup \left\{ d(f(i),g(i)) \mid (i,d) \in A \right\}.
\]
Each projection $\pi_{i} \colon \prod_{i \in I}X_i \to X_{i}$ is
obviously a homomorphism.
We say that the product of a family $(X_{i})_{i \in I}$ is
\emph{inhabited} if the underlying set $\prod_{i \in I}X_i$ is
inhabited.
\begin{prop}\label{prop:OmegaProd}
  The functor $\overline{(-)} \colon \GUS \to \FTop$ preserves
  inhabited countable products.
\end{prop}
\proof
Let $\left( (X_{n},M_{n})\right)_{n \in \Nat}$ be a sequence of gus's
with a chosen sequence $\varphi \in \prod_{n \in \Nat}X_{n}$.
Let $\prod_{n \in \Nat}X_{n}$ denote the product of the family,
where we left the underlying family of generalised metrics implicit.
Write $\mathcal{U}(\prod_{n \in \Nat}X_{n}) = (U_{X},
\cov_{X}, \leq_{X})$ for its localic completion.
The elements of $\prod_{n \in \Nat}X_{n}$ will be denoted by Greek
letters $\alpha,\beta,\gamma$, and we write $\alpha_{n}$ for
$\alpha(n)$.

Given a sequence $\left(r_{n} \colon \mathcal{S} \to
\mathcal{U}(X_{n})  \right)_{n \in \Nat}$ of formal topology maps, 
defined  a relation $r \subseteq S \times U_{X}$ by
\[
  a \mathrel{r} \ball_{A}(\alpha,\varepsilon)
  \defeqiv \left( \exists \ball_{B}(\beta,\delta) <_{X}
  \ball_{A}(\alpha,\varepsilon)\right) \left( \forall (i, d) \in B \right)
  a \cov r_{i}^{-}\ball_{d}(\beta_{i},\delta).
\]
We claim that $r$ is a formal topology map. We only show that $r$ satisfies
\ref{FTM2} since other conditions are easy to check. 
Let $\ball_{A}(\alpha,\varepsilon), \ball_{B}(\beta,\delta) \in U_{X}$,
and let
$a \in r^{-}\ball_{A}(\alpha,\varepsilon) \downarrow
r^{-}\ball_{B}(\beta,\delta)$. Then, there exist
$\ball_{A'}(\alpha',\varepsilon') <_{X} \ball_{A}(\alpha,\varepsilon)$
and $\ball_{B'}(\beta', \delta') <_{X} \ball_{B}(\beta,\delta)$ such
that 
\begin{itemize}
  \item $\left( \forall (i,d) \in A' \right) a \cov
    r_{i}^{-}\ball_{d}(\alpha_{i}', \varepsilon')$,

  \item $\left( \forall (j,\rho) \in B' \right) a \cov
    r_{j}^{-}\ball_{\rho}(\beta_{j}', \delta')$.
\end{itemize}
We can write $A'$ as a disjoint union
\begin{equation*}
  A' = \left( \left\{ i_{0} \right\} \times A_{0} \right) \cup
  \cdots \cup \left( \left\{ i_{n} \right\} \times A_{n} \right),
\end{equation*}
where for each $k \leq n$, we have $i_{k} \in \Nat$ and $A_{k} \in
\PFin(M_{i_{k}})$, and $0 \leq k < k' \leq n \implies i_{k} \neq i_{k'}$.
Similarly, write $B'$ as a disjoint union
  $
  B' = \left(\left\{ j_{0} \right\} \times B_{0}  \right) \cup
  \cdots \cup \left( \left\{ j_{m} \right\} \times B_{m} \right)
  $
that satisfies the analogous properties as those of $A'$.
Put $I = \left\{ i_{k} \mid k \leq n \right\}$,
$J = \left\{ j_{k} \mid k \leq m \right\}$,
and $P = I \cup J = \left\{ p_{0}, \dots, p_{l} \right\}$ with the
property $0 \leq k < k' \leq l \implies p_{k} \neq p_{k'}$.

Choose $\theta \in \PRat$ such that $A(\alpha,\alpha') + \varepsilon'
+ 2 \theta < \varepsilon$, and 
$B(\beta, \beta') + \delta' + 2 \theta < \delta$.
For each $k \leq l$, define a subset $V_{k} \subseteq
U_{X_{p_{k}}}$ by cases:
\begin{enumerate}
  \item If $p_{k} \in I \cap J$, then put
      $
      V_{k} \defeql V(A_{p_{k}}) \downarrow V(B_{p_{k}})
      \downarrow
      \mathcal{C}^{\theta}_{\sup (A_{p_{k}} \cup B_{p_{k}})}.
      $
    Here, the subset $V(A_{p_{k}}) \subseteq U_{p_{k}}$
    is defined by 
    \[
      V(A_{p_{k}}) \defeql
      \ball_{d_{0}}(\alpha_{p_{k}}',\varepsilon')
      \downarrow \cdots \downarrow
      \ball_{d_{N}}(\alpha_{p_{k}}',\varepsilon'),
    \]
    where $A_{p_{k}} = \left\{ d_{0},\dots,d_{N} \right\}$.
    The subset $V(B_{p_{k}}) \subseteq U_{p_{k}}$ is defined similarly.

  \item If $p_{k} \in I \setminus J$,  put
      $
      V_{k} \defeql V(A_{p_{k}}) \downarrow
      \mathcal{C}^{\theta}_{\sup (A_{p_{k}})}.
      $

  \item If $p_{k} \in J \setminus I$, put
      $
      V_{k} \defeql V(B_{p_{k}}) \downarrow
      \mathcal{C}^{\theta}_{\sup (B_{p_{k}})}.
      $
\end{enumerate}
Then, we have $a \cov r_{p_{0}}^{-} V_{0} \downarrow \cdots \downarrow
r_{p_{l}}^{-} V_{l}$. Let $b \in r_{p_{0}}^{-} V_{0} \downarrow
\cdots \downarrow r_{p_{l}}^{-} V_{l}$. Then, for each $k \leq l$, there
exists $\ball_{\sigma_{k}}(z_{k}, \xi_{k}) \in V_{k}$ such that
$b \cov r_{p_{k}}^{-}\ball_{\sigma_{k}}(z_{k},\xi_{k})$. Define
$\gamma \in \prod_{n \in \Nat} X_{n}$ by
\[
  \gamma_{n} = 
  \begin{cases}
    z_{k}& \text{$n = p_{k}$ for some $k \leq l$},\\
    \varphi_{n}& \text{otherwise},
  \end{cases}
\]
and put $C = \left\{ (p_{k},\sigma_{k}) \mid k \leq l\right\}$.
Then, we have $b \mathrel{r} \ball_{C}(\gamma, 2 \theta)$. Moreover,
\begin{align*}
  A(\alpha,\gamma) + 2 \theta
  &\leq A(\alpha,\alpha') + A(\alpha',\gamma) + 2 \theta\\
  &\leq A(\alpha,\alpha') + A'(\alpha',\gamma) + 2 \theta\\
  &<  A(\alpha,\alpha') + \varepsilon' + 2 \theta < \varepsilon.
\end{align*}
Thus $\ball_{C}(\gamma,2\theta) <_{X} \ball_{A}(\alpha,\varepsilon)$.
Similarly we have $\ball_{C}(\gamma,2\theta) <_{X}
\ball_{B}(\beta,\delta)$. Hence, 
\[
  a \cov r^{-}(\ball_{A}(\alpha,\varepsilon) \downarrow
  \ball_{B}(\beta,\delta)).
\]

Next, we note that the following holds:
\[
  a \mathrel{r} \ball_{A}(\alpha,\varepsilon) \implies \left( \forall
  (i,d) \in A \right) a \cov r_{i}^{-}\ball_{d}(\alpha_{i},
  \varepsilon).
\]
The proof is similar to the above proof of the condition \ref{FTM2}
for $r$. Moreover, the functor $\overline{(-)} \colon \GUS \to
\FTop$ sends each projection $\pi_{n} \colon \prod_{n \in
\Nat}X_{n} \to X_{n}$ to a formal topology map
$\overline{\pi_{n}} \colon \mathcal{U}(\prod_{n \in \Nat}X_{n}) \to
\mathcal{U}(X_{n})$. Clearly, we have the following:
\[
  \ball_{A}(\alpha,\varepsilon) \mathrel{\overline{\pi_{n}}}
  \ball_{d}(x,\delta) \iff \ball_{A}(\alpha,\varepsilon) <_{X}
  \ball_{ \left\{(n,d) \right\}}(\varphi_{(n,x)},\delta),
\]
where $\varphi_{(n,x)}$ is obtained from $\varphi$ by replacing
$n$th element with $x$.

With these at our disposal, it is straightforward to show that
$r$ makes the diagram
\begin{equation*}
  \xymatrix{
    \mathcal{S}
    \ar[rd]_{r_{n}} \ar[r]^-{r} &
    \mathcal{U}(\prod_{n \in \Nat}X_n )
    \ar[d]^{\overline{\pi_{n}}} &\\
    & \mathcal{U}(X_n)
  }
\end{equation*}
commute, and it is a unique such morphism.
\qed
Note that Proposition \ref{prop:OmegaProd} generalises the fact that
the localic completion of Baire space $\Nat^{\Nat}$ is the point-free
Baire space \cite[Proposition 3.1]{PalmgrenFormalContUnifContNCI}.

Let $(X,M)$ be a gus, and let $M^{\op}$ denote $M$ ordered
by the opposite of \eqref{eq:OrderGM}.  Then, we have a
cofiltered diagram $D_{M} \colon M^{\op} \to \GUS$ given by $D_{M}(d)  = (X,
d)$ for $d \in M$, and for each $d,d' \in M$ such that $d
\leq d'$ we have $D_{M}(d,d') \defeql \id_{X} \colon (X,d') \to 
(X,d)$.
For each $d \in M$, we have a homomorphism $\sigma_{d} \colon (X,M)
\to (X,d)$ with the identity function $\id_{X}$ as the underlying map.  It is
easy to see that the family $\left( \sigma_d \colon (X,M) \to
(X,d)\right)_{d \in M}$ is a limit of the diagram  $D_{M}$.

\begin{prop}\label{prop:PresColim}
  For any gus $(X,M)$, the functor
  $\overline{(-)} \colon \GUS \to \FTop$ preserves
  the limit of the cofiltered diagram $D_{M} \colon M^{\op} \to \GUS$.
\end{prop}
\proof
We write $\mathcal{U}(X)$ for the localic completion of $(X,M)$
and $\mathcal{M}(X,d)$ for the localic completion of $(X,d)$ for each
$d \in M$.
It suffices to show that the family 
\[
  \left( \overline{\sigma_d} \colon
  \mathcal{U}(X) \to \mathcal{M}(X,d)\right)_{d \in M}
\]
is a limit of
the diagram $\overline{(-)} \circ D_{M}  \colon M^{\op} \to \FTop$.

First, for each $d,d' \in M$ such that $d \leq d'$, we have
\[
  \ball_{d'}(x',\varepsilon') \mathrel{\overline{D_{M}(d,d')}}
  \ball_{d}(x,\varepsilon) \iff \ball_{d'}(x',\varepsilon')
  <_{X} \ball_{d}(x,\varepsilon).
\]
Moreover, for each $d \in M$, we have
\[
  \ball_{d'}(x',\varepsilon') \mathrel{\overline{\sigma_{d}}}
  \ball_{d}(x,\varepsilon) \iff \ball_{d'}(x',\varepsilon')
  <_{X} \ball_{d}(x,\varepsilon).
\]
Given any cone $\left( r_{d} \colon \mathcal{S} \to
\mathcal{M}(X,d)\right)_{d \in M}$ over $\overline{(-)} \circ
D_{M}  \colon M^{\op} \to \FTop$, define a relation
$r \subseteq S \times U_X$ by
\[
  a \mathrel{r} \ball_{d}(x,\varepsilon) \defeqiv a
  \mathrel{r_{d}} \ball_{d}(x,\varepsilon)
\]
for all $a \in S$ and $\ball_{d}(x,\varepsilon) \in U_X$.
It is straightforward to show that $r$ is a formal topology map
$r \colon \mathcal{S} \to \mathcal{U}(X)$. For example, the property
\ref{FTM1} follows from the fact that $M$ is inhabited.
For the property \ref{FTM2}, let
$\ball_{d_{1}}(x_{1},\varepsilon_{1}),\ball_{d_{2}}(x_{2},\varepsilon_{2})
\in U_{X}$.  Putting $d = \sup\left\{ d_{1}, d_{2} \right\}$,  we have
\begin{align*}
  r^{-} \ball_{d_{1}}(x_{1},\varepsilon_{1}) \downarrow r^{-}
\ball_{d_{2}}(x_{2},\varepsilon_{2})
&\cov 
r_{d}^{-}\overline{D_{M}(d_{1},d)}^{-}
\ball_{d_{1}}(x_{1},\varepsilon_{1}) \downarrow
r_{d}^{-}\overline{D_{M}(d_{2},d)}^{-}
\ball_{d_{2}}(x_{2},\varepsilon_{2}) \\
&\cov 
r_{d}^{-}\left(\overline{D_{M}(d_{1},d)}^{-}
\ball_{d_{1}}(x_{1},\varepsilon_{1}) \downarrow
\overline{D_{M}(d_{2},d)}^{-}
\ball_{d_{2}}(x_{2},\varepsilon_{2}) \right) \\
&\cov 
r^{-}\left( \ball_{d_{1}}(x_{1},\varepsilon_{1}) \downarrow 
\ball_{d_{2}}(x_{2},\varepsilon_{2}) \right).
\end{align*}
The other properties of $r$ are easy to check. Then, it is
straightforward to show that $ \overline{\sigma_{d}}\circ r = r_{d}$
for each $d \in M$, and that $r$ is the unique formal topology map
with this property.
\qed

\subsection{Symmetric generalised uniform spaces}\label{sec:SymmGUS}
We fix a symmetric gus $(X,M)$ throughout this subsection.
The aims of this subsection are twofold. The first is to obtained
a point-free analogue of the point-set completion of $(X,M)$
as a closed subspace of the product of the completion of $(X,d)$
for each $d \in M$ (Section \ref{sec:Closed Embedding}). The second is to analyse the point-free uniform
structure on $\mathcal{U}(X)$ induced by each symmetric generalised
metric $d \in M$ (Section \ref{sec:PFUniStrct}), and relate it to the complete uniform structure on
$\Pt(\mathcal{U}(X))$ (Section \ref{sec:PtLComp}).

\subsubsection{Closed embedding into $\prod_{d \in
M}\mathcal{M}(X,d)$}\label{sec:Closed Embedding}
Since $\mathcal{U}(X)$ is a limit of the diagram
$\overline{(-)} \circ  D_{M}  \colon M^{\op} \to \FTop$,
the unique morphism $\iota_{\mathcal{U}(X)} \colon \mathcal{U}(X)
\to \prod_{d \in M}\mathcal{M}(X,d)$ determined by 
the family
$\left( \overline{\sigma_{d}} \colon \mathcal{U}(X) \to
\mathcal{M}(X,d) \right)_{d \in M}$ is an embedding (because
it must be an equaliser). The image of 
$\iota_{\mathcal{U}(X)}$ is overt with positivity
\[
  \iota_{\mathcal{U}(X)}[U_{X}] = \left\{ A \in \Fin\Bigl(\sum_{d \in
  M}U_{(X,d)}\Bigr) \,{\Big\vert}\, \left( \exists a \in U_{X} \right) \left( \forall
  (d, b) \in A \right) a \cov_{X} b \right\}.
\]

\begin{lem}
  $
  \iota_{\mathcal{U}(X)}[U_{X}]
  =
  \left\{ A \in \Fin\bigl(\sum_{d \in M}U_{(X,d)}\bigr) 
  \mid \left( \exists x  \in X \right) \left( \forall
  (d, a) \in A \right) x \in a_{*} \right\}.
  $
\end{lem}
\proof
The inclusion $\subseteq$ is clear. Conversely, let 
$A \in \Fin\bigl(\sum_{d \in M}U_{(X,d)}\bigr)$, and let $x \in X$
such that  $x \in a_{*}$ for each $(d, a) \in A$. Without loss of
generality, we can assume that $A$ is inhabited. Write 
$A = \left\{ (d_{0},\ball_{d_{0}}(x_{0},\varepsilon_{0})),\dots,(d_{n},\ball_{d_{n}}(x_{n},\varepsilon_{n}))
\right\}$, and choose $\theta \in \PRat$ such that $d_{i}(x_{i}, x) +
\theta < \varepsilon_{i}$ for each $i \leq n$.
Let $d = \sup \left\{ d_{i} \mid i \leq n \right\}$.  Then, $\ball_{d}(x,\theta)
  <_{X} \ball_{d_{i}}(x_{i},\varepsilon_{i})$ for each $i \leq n$. Hence
  $A \in \iota_{\mathcal{U}(X)}[U_{X}]$.
\qed
Define $W \subseteq \Fin\bigl(\sum_{d \in M}U_{(X,d)}\bigr)$  by
  $
    W \defeql  \iota_{\mathcal{U}(X)}[U_{X}].
  $
Let $\mathcal{S}_{W} = (S_{W}, \cov_{W} \leq)$ be the overt weakly
closed subtopology of $\prod_{d \in M}\mathcal{M}(X,d)$ determined by
$W$, where $\leq$ is the preorder on the base of $\prod_{d \in M}\mathcal{M}(X,d)$.

\begin{lem}
  For each $d,d' \in M$ such that $d \leq d'$, the diagram commutes:
  \[
  \xymatrix{
    & \mathcal{S}_{W} \ar[ld]_{p_{d'}} \ar[rd]^{p_{d}}&\\
    \mathcal{M}(X,d') \ar[rr]^{\overline{D_{M}(d,d')}} & &
    \mathcal{M}(X,d)
      }
  \]
  where $p_{d}$ (respectively $p_{d'}$) is the restriction of
  the projections $p_{d} \colon \prod_{d \in M}\mathcal{M}(X,d) \to
  \mathcal{M}(X,d)$ to $\mathcal{S}_{W}$.
\end{lem}
\proof
Let $\ball_{d}(x,\varepsilon) \in U_{(X,d)}$ and
$A \in S_{W}$.

First, suppose that $A \mathrel{p_{d}} \ball_{d}(x,\varepsilon)$, i.e.\
$A = \left\{ (d, \ball_{d}(x,\varepsilon)) \right\}$.
We must show that 
\[
  A \cov_{W} \left\{ \left\{ (d', \ball_{d'}(y,\delta)) \right\}  \mid
  \ball_{d'}(y,\delta) <_{X} \ball_{d}(x,\varepsilon) \right\}.
\]
By \ref{U1}, we have
\[
  A \cov_{W} \left\{ \left\{(d, \ball_{d}(x,\varepsilon')) \right\} \mid \varepsilon'
< \varepsilon \right\}.
\]
Let $\varepsilon' \in \PRat$ such that $\varepsilon' < \varepsilon$,
and choose $\theta \in \PRat$ such that $\varepsilon' + 2\theta <
\varepsilon$. By \ref{U2}, we have
\[
  \left\{ (d,\ball_{d}(x,\varepsilon')) \right\} \cov_{W}
  \left\{ \left\{ (d,\ball_{d}(x,\varepsilon')),
  (d',\ball_{d'}(y,\theta))  \right\} \mid y \in X \right\} \cap W.
\]
Let $y \in X$, and suppose that 
$\left\{ (d,\ball_{d}(x,\varepsilon')),
  (d',\ball_{d'}(y,\theta))  \right\} \in W$. Then, there exists
  $z \in X$ such that $d'(y,z) < \theta$ and $d(x,z) < \varepsilon'$.
  Thus, $d(x,y) + \theta \leq
  d(x,z) + d(z,y) + \theta < 
  d(x,z) + d'(y,z) + \theta < \varepsilon' + 2\theta < \varepsilon$. Hence
  $ \left\{ (d,\ball_{d}(x,\varepsilon')),
  (d',\ball_{d'}(y,\theta)) \right\} \leq \left\{ (d',
  \ball_{d'}(y, \theta))\right\}$ and $\ball_{d'}(y,\theta)
  <_{X} \ball_{d}(x,\varepsilon)$.
 
  The proof of the converse, i.e.\ $A \mathrel{(\overline{D_{M}(d,d')} \circ
  p_{d'})} \ball_{d}(x,\varepsilon) \implies A \cov_{W}
  p_{d}^{-}\ball_{d}(x,\varepsilon)$ is similar.
\qed

By Proposition \ref{prop:PresColim}, there exists a unique morphism $s \colon \mathcal{S}_{W} \to
\mathcal{U}(X)$ such that $p_{d} = \overline{\sigma_{d}} \circ s$
for each $d \in M$. By the proof of Proposition \ref{prop:PresColim},
we have
\[
  A \mathrel{s} \ball_{d}(x,\varepsilon) \iff A \cov_{W} \left\{
    (d, \ball_{d}(x,\varepsilon)) \right\}.
\]
It is straightforward to show that $\iota_{\mathcal{U}(X)} \circ s =
\id_{\mathcal{S}_{W}}$. Since $\iota_{\mathcal{U}(X)}$ is a
monomorphism, it is an isomorphism.
\begin{thm}\label{thm:ClosedEmbedding}
  The localic completion $\mathcal{U}(X)$ embeds into
  $\prod_{d \in M} \mathcal{M}(X,d)$ as an overt weakly closed
  subtopology.
\end{thm}

\subsubsection{Point-free uniform structures}\label{sec:PFUniStrct}
Vickers \cite[Proposition 6.6]{LocalicCompletionGenMet} showed that
each symmetric generalised metric $d \in M$ determines a morphism
$\overline{d} \colon \mathcal{M}(X,d) \times \mathcal{M}(X,d) \to \FUR$ of formal topologies
determined  by a function $f_{d} \colon \PRat \to \Pow(U_{(X,d)} \times
U_{(X,d)})$ given by
\[
  f_{d}(q) \defeql \left\{  
  \left(
  \ball_{d}(x_{_1},\varepsilon_{1}),\ball_{d}(x_{_2},\varepsilon_{2})
  \right) \mid 
  d(x_{1},x_{2}) + \varepsilon_{1} + \varepsilon_{2} < q \right\}.
\]
He showed that $\overline{d}$ satisfies the properties of point-free 
symmetric generalised metrics.\footnote{Vickers's proof is in the
setting of impredicative topos theory, but it is straightforward to
adapt his proof to the setting of formal topology.} This means that
\begin{enumerate}
  \item  $\overline{d} \circ \Delta_{\mathcal{M}(X,d)} = \mathbf{0}$,
  \item  $\overline{d} \circ \tau_{\mathcal{M}(X,d)} = \overline{d}$,
  \item  $\langle \overline{d} \circ \langle p_{1} \circ q_{1},
    p_{2} \circ q_{2} \rangle,  + \circ \langle \overline{d} \circ q_{1},
    \overline{d} \circ q_{2} \rangle \rangle$
    factors uniquely through the embedding $\leq_{u} \to \FUR
    \times \FUR$.
\end{enumerate}
Here 
\begin{itemize}
  \item $\Delta_{\mathcal{M}(X,d)}$ is the diagonal morphism
        $\langle \id_{\mathcal{M}(X,d)}, \id_{\mathcal{M}(X,d)} \rangle \colon
        \mathcal{M}(X,d) \to \mathcal{M}(X,d) \times \mathcal{M}(X,d)$;

  \item
  $\mathbf{0} \colon \mathcal{M}(X,d) \to \FUR$ is determined by
  a function $f_{\mathbf{0}} \colon \PRat \to \Pow(U_{(X,d)})$  given by 
    $
    f_{\mathbf{0}}(q) = U_{(X,d)}
    $ for all $q \in \PRat$;

  \item $\tau_{\mathcal{M}(X,d)}$ is
    the twisting morphism $\langle p_{2}, p_{1} \rangle \colon
    \mathcal{M}(X,d) \times \mathcal{M}(X,d) \to \mathcal{M}(X,d) \times \mathcal{M}(X,d)$, where
    $p_{1}$ and $p_{2}$ are the first and second projections of
    $\mathcal{M}(X,d) \times \mathcal{M}(X,d)$;

  \item $q_{1}, q_{2} \colon \mathcal{T} \to \mathcal{M}(X,d) \times
    \mathcal{M}(X,d)$ are the pullback:
  \[
    \xymatrix{
      \mathcal{T} \ar[r]^-{q_{2}} \ar[d]_-{q_{1}} &\mathcal{M}(X,d)
      \times \mathcal{M}(X,d) \ar[d]^-{p_{1}} \\
     \mathcal{M}(X,d)  \times \mathcal{M}(X,d) \ar[r]^-{p_{2}} & \mathcal{M}(X,d) 
      }
  \]
\end{itemize}
Moreover, Vickers showed that if $d$ is Dedekind,
then $\overline{d}$ uniquely factors through the Dedekind reals
$\FNNR$ via $\iota_D \colon \FNNR \to \FUR$ (his result can be
easily adapted to the finite Dedekind reals). If $\overline{d}_{D} \colon
\mathcal{M}(X,d) \times \mathcal{M}(X,d) \to \FNNR$ is a factorisation
of $\overline{d}$, then by using the facts in Section
\ref{sec:UpperReals}, it is easy to see that $\overline{d}_{D}$
satisfies the properties of point-free (non-finite) Dedekind symmetric
generalised metrics on  $\mathcal{M}(X,d)$.

Furthermore, for any formal topology map $r \colon \mathcal{S} \to
\mathcal{M}(X,d)$, a straightforward diagram chasing shows that  $\overline{d} \circ
r \times r$ (or $\overline{d}_{D} \circ r \times r$) satisfies the properties of
point-free (non-finite Dedekind) symmetric generalised metric on $\mathcal{S}$.
Thus, each symmetric generalised metric $d \in M$ gives rise to a
point-free (non-finite Dedekind) symmetric generalised metric on
$\mathcal{U}(X)$ by composing $\overline{d} \colon \mathcal{M}(X,d)
\times \mathcal{M}(X,d) \to \FUR$ (or $\overline{d}_{D}$) with the
products of projection $p_{d} \colon \prod_{d \in M}
\mathcal{M}(X,d) \to \mathcal{M}(X,d)$ and the embedding
$\iota_{\mathcal{U}(X)} \colon \mathcal{U}(X) \to \prod_{d \in M}
\mathcal{M}(X,d)$.
Let us denote this composite by
\[
  \widetilde{d} \defeql \overline{d} \circ p_{d} \times p_{d} \circ
  \iota_{\mathcal{U}(X)} \times \iota_{\mathcal{U}(X)}.
\]
It is easy to see that $\widetilde{d}$ is determined by a function
$g_{d} \colon \PRat \to \Pow(U_{X} \times U_{X})$ given by
\[
  g_{d}(q) \defeql \left\{ \left(\ball_{d}(x_{1},\varepsilon_{1}),
  \ball_{d}(x_{2},\varepsilon_{2})\right) \in U_X \times U_X  \mid
  d(x_{1},x_{2}) +
  \varepsilon_{1} + \varepsilon_{2} < q \right\}.
\]
Moreover, if we regard $(X,M)$ as a discrete formal topology (with
base $X$ with the discrete order $=$ and with the trivial cover $x \cov U \iff x \in U$), then each $d \in M$ determines a
point-free symmetric generalised metric $d \colon X \times X \to \FUR$
on $X$, which is given by a function $h_{d} \colon \PRat \to \Pow(X \times
X)$ defined by
\[
  h_{d}(q) \defeql \left\{ (x,x') \in X \times X  \mid d(x,x') < q \right\}.
\]
Also, there is a formal topology map $\iota_{X} \colon X \to \mathcal{U}(X)$
defined by 
\[
  x \mathrel{\iota_{X}} \ball_{d}(x',\varepsilon) \defeqiv d(x',x) <
  \varepsilon,
\]
which is easily seen to be an isometry in the sense that $d =
\widetilde{d} \circ \iota_X \times \iota_X$ for each $d \in M$.
  
\subsubsection{Formal points of localic completions}\label{sec:PtLComp}
In this subsection, we work impredicatively, assuming that
$\Pow(X)$ is a set for any set $X$. 
When specialised to finite Dedekind symmetric gus's, however, 
the results in this subsection give a predicative completion of a given
gus (see Remark \ref{rem:Predicativity}).
In the next subsection (Section \ref{sec:Predicative}), we fully
address the predicativity issue raised in this subsection.
As in Section \ref{sec:PFUniStrct}, we work on a fixed symmetric
generalised uniform space $(X,M)$.

For each $d \in M$, by applying the operation $\Pt(-)$ to the
point-free symmetric generalised metric $\widetilde{d} \colon
\mathcal{U}(X) \times \mathcal{U}(X) \to \FUR$, we obtain a
symmetric generalised metric 
\[
  \Pt(\widetilde{d}) \colon
  \Pt(\mathcal{U}(X)) \times \Pt(\mathcal{U}(X)) \to \UR
\]
on $\Pt(\mathcal{U}(X))$. This is because $\Pt(-)$ is a right adjoint
and so it preserves all the properties of symmetric generalised
metrics. By an abuse of notation, we write $\widetilde{d}$ for
$\Pt(\widetilde{d})$.
Explicitly, $\widetilde{d} \colon \Pt(\mathcal{U}(X)) \times
\Pt(\mathcal{U}(X)) \to \UR$ is given by
\begin{align*}
  \widetilde{d}(\alpha,\beta) = \left\{ q \in \PRat \mid \left(
  \exists \ball_{d}(x,\varepsilon) \in \alpha\right)\left( \exists
  \ball_{d}(y,\delta) \in \beta\right) d(x,y) + \varepsilon +
  \delta < q \right\}.
\end{align*}
Moreover, $\widetilde{d}$ factors
through (finite) Dedekind reals if $d$ is (finite) Dedekind.
Note that $d \leq \rho \iff \widetilde{d} \leq \widetilde{\rho}$
for all $d, \rho \in M$, and $\widetilde{\sup\left\{ d,\rho \right\}}
= \sup\left\{ \widetilde{d}, \widetilde{\rho} \right\}$.

On the other hand, under the operation $\Pt(-)$, the gus $(X,M)$ as a
discrete formal
topology (with the formal uniform structure $M$) is mapped essentially
to $(X,M)$ itself. The formal topology map $\iota_X \colon X \to
\mathcal{U}(X)$ is mapped by $\Pt(-)$ to a function
$i_X \colon X \to \Pt(\mathcal{U}(X))$ given by
  \begin{equation}
    \label{eq:Completion}
  i_X(x) \defeql \Diamond x.
  \end{equation}
  Moreover, it is an isometry, i.e. $d = \widetilde{d} \circ
  i_{X} \times i_{X}$ for each $d \in M$.

\begin{lem}\label{lem:MiscDense}
  For each $\alpha \in \Pt(\mathcal{U}(X))$, we have
  \[
    \ball_{d}(x,\varepsilon) \in \alpha \iff \widetilde{d}(\alpha,
    \Diamond x) < \varepsilon.
  \]
\end{lem}
\proof
Suppose that $\ball_{d}(x,\varepsilon) \in \alpha$.
By \ref{U1}, there exists $\varepsilon' < \varepsilon$
such that $\ball_{d}(x,\varepsilon') \in \alpha$.
Choose $\theta \in \PRat$ such that $\varepsilon' + \theta <
\varepsilon$. 
Since $\ball_{d}(x,\theta) \in \Diamond x$, we have
$\widetilde{d}(\alpha, \Diamond x) < \varepsilon$.
Conversely, suppose that $\widetilde{d}\left(\alpha, \Diamond x \right)
< \varepsilon$. Then, there exist $\ball_{d}(y,\delta) \in \Diamond x$
and $\ball_{d}(z,\xi) \in \alpha$ such that $d(y,z) + \delta +
\xi < \varepsilon$. Then, $d(x,z) + \xi \leq d(x,y) + d(y,z) +
\xi < \delta + d(y,z) + \xi < \varepsilon$, and hence
$\ball_{d}(z,\xi) <_{X} \ball_{d}(x,\varepsilon)$. Therefore,
$\ball_{d}(x,\varepsilon) \in \alpha$.
\qed

\begin{prop}\label{prop:DenseEmb}
  The function $i_{X} \colon X \to \Pt(\mathcal{U}(X))$ is dense, i.e.
  \[
    \left( \forall \alpha \in \Pt(\mathcal{U}(X)) \right) \left(
    \forall (d,\varepsilon) \in \Ent(M) \right)\left( \exists x
    \in X \right) \widetilde{d}(\alpha, \Diamond x) < \varepsilon.
  \]
\end{prop}
\proof
By the axiom \ref{U2} and Lemma \ref{lem:MiscDense}.
\qed

\begin{defi}
  A symmetric gus $(X,M)$ is \emph{separated} if 
  \[
    \left[ \left( \forall d \in M \right) d(x,y) = 0 \right] \iff x = y.
\]
\end{defi}
\begin{prop}\label{prop:SeparatedGUS}
  A symmetric gus $(X,M)$ is separated if and only if the function  $i_{X} \colon X \to
  \Pt(\mathcal{U}(X))$ is injective.
\end{prop}
\proof
Suppose that $(X,M)$ is separated, and assume $i_{X}(x) =
i_{X}(y)$. Let $d \in M$. Then for each $q \in \PRat$, since
$\ball_{d}(x,q) \in i_{X}(x)$, we have
$d(x,y) < q$. Thus, $d(x,y) = 0$. Hence $x = y$.

Conversely, suppose that $i_{X}$ is injective. Let $x,y \in X$, and
suppose that $d(x,y) = 0$ for all $d \in M$. Let
$\ball_{d}(z,\varepsilon) \in i_{X}(x)$. Then, $d(z,x) <
\varepsilon$, and so $d(z,y) < \varepsilon$. Thus,
$\ball_{d}(z,\varepsilon) \in i_{X}(y)$. Similarly (using symmetry),
we have $i_{X}(y) \subseteq i_{X}(x)$. Since $i_{X}$ is injective,
we have $x = y$.
\qed

In the rest of this subsection, we show that $\Pt(\mathcal{U}(X))$
equipped with $\widetilde{M} \defeql \left\{
  \widetilde{d} \mid d \in M \right\}$ is a completion of $(X,M)$.

\begin{defi}\label{def:CauchyFilterUSpa}
  A \emph{Cauchy filter} on a symmetric gus $(X,M)$ is a
  set $\mathcal{F}$ of subsets of $X$ such that 
  \begin{enumerate}[label=(CF\arabic*)]
    \item\label{CF1} $ U \in \mathcal{F}  \implies U \meets X$,
    \item\label{CF2} $ U \in \mathcal{F} \amp U \subsets V \implies V \in
      \mathcal{F}$,
    \item\label{CF3} $ U,V \in \mathcal{F} \implies
      U \cap V \in \mathcal{F}$,
    \item\label{CF4} $\left( \forall \left(d,\varepsilon \right) \in
      \Ent(M)
    \right)\left( \exists x \in X \right)  \Ball_d(x,\varepsilon) \in
    \mathcal{F}$.
  \end{enumerate}
  A Cauchy filter $\mathcal{F}$ on $X$ \emph{converges} to a point
  $x \in X$ if
    $
    \left(\forall (d,\varepsilon) \in \Ent(M) \right)
    \Ball_{d}(x,\varepsilon) \in \mathcal{F}.
    $
  A symmetric gus is \emph{complete} if every Cauchy filter
  on $X$ converges to some point.
\end{defi}

\begin{lem}\label{lem:FiltPt}
  Let $(X,M)$ be a symmetric gus.
  \begin{enumerate}
    \item \label{lem:FiltPt1} If $\mathcal{F}$ is a Cauchy filter on $X$, then
      \[
        \alpha_{\mathcal{F}} \defeql \left\{ a \in U_X \mid \left( \exists b
          <_{X} a \right) b_{*} \in \mathcal{F} \right\}
        \]
        is a formal point of $\mathcal{U}(X)$.
      \item \label{lem:FiltPt2} If $\alpha$ is a formal point of $\mathcal{U}(X)$, then
        \[
          \mathcal{F}_{\alpha}  \defeql \left\{ U \in \Pow(X) \mid \left(
            \exists  a \in \alpha \right) a_{*} \subseteq U \right\}
          \]
          is a Cauchy filter on $X$. Moreover,
          $\alpha_{\mathcal{F}_{\alpha}} = \alpha$.
  \end{enumerate}
\end{lem}
\proof
\eqref{lem:FiltPt1} For example, to see that
$\alpha_{\mathcal{F}}$ satisfies \ref{P2}, let $a, b \in
\alpha_{\mathcal{F}}$. Then, there exist $a' <_{X} a$ and
$b'<_{X} b$ such that $a'_{*}, b'_{*} \in \mathcal{F}$.
Write
$a = \ball_{d_1}(x,\varepsilon)$,
$b = \ball_{d_2}(y,\delta)$,
$a' = \ball_{d_1'}(x',\varepsilon')$, and
$b' = \ball_{d_2'}(y',\delta')$. Choose $\theta \in \PRat$ such that
$\ball_{d_1'}(x',\varepsilon' + 3\theta) <_{X} a$, and
$\ball_{d_2'}(y',\delta' + 3\theta) <_{X} b$, and put $\rho = \sup
\left\{ d_1', d_2' \right\}$. Since
$\mathcal{F}$ is Cauchy, there exists $z \in X$ such that 
$\Ball_{\rho}(z,\theta) \in \mathcal{F}$. Thus,
$\Ball_{\rho}(z,\theta) \meets \Ball_{d_1'}(x',\varepsilon')$ and 
$\Ball_{\rho}(z,\theta) \meets \Ball_{d_2'}(y',\delta')$.
Then, 
\[
  d_1(x,z) + 2\theta \leq d_1(x,x') + d_1(x',z) + 2\theta < d_1(x,x') +
  \varepsilon' + \theta + 2\theta < \varepsilon.
\]
Thus, $\ball_{\rho}(z, 2\theta) <_{X} a$, and similarly
$\ball_{\rho}(z, 2\theta) <_{X} b$. Moreover,  $\ball_{\rho}(z,
2\theta) \in \alpha_{\mathcal{F}}$.

\eqref{lem:FiltPt2}
The claim that  $ \mathcal{F}_{\alpha}$ is a Cauchy filter is obvious.
To see that $\alpha_{\mathcal{F}_{\alpha}} = \alpha$, let
$a = \ball_{d}(x,\varepsilon) \in \alpha_{\mathcal{F}_{\alpha}}$.
Then, there exists $\ball_{d'}(x',\varepsilon') <_{X} a$ such that $\Ball_{d'}(x',\varepsilon') \in
\mathcal{F}_{\alpha}$. Thus, there exists $\ball_{\rho}(y,\delta) \in \alpha$
such that $\Ball_{\rho}(y, \delta) \subseteq \Ball_{d'}(x',\varepsilon')$.
Choose $\theta \in \PRat$ such that $d(x,x') + \varepsilon' + \theta <
\varepsilon$. Then, there exists $c = \ball_{\rho'}(y',\delta') \in
\alpha$ such  that $c \in
\ball_{\rho}(y,\delta) \downarrow \mathcal{C}^{\theta}_{d'}$.
Then
\[
  d(x,y') + \delta' \leq d(x,x') + d(x',y') + \delta'
  < d(x,x') + \varepsilon' + \theta < \varepsilon.
\]
Thus, $c <_{X} a$, and hence $a \in \alpha$.

Conversely, if $a \in \alpha$, then there exists $b <_{X} a$ such
that $b \in \alpha$ by \ref{U1}. Then, $b_{*} \in
\mathcal{F}_{\alpha}$ and so $a \in \alpha_{\mathcal{F}_{\alpha}}$.
\qed

\begin{lem}\label{lem:FiltConv}
  A Cauchy filter $\mathcal{F}$ on $(X,M)$ converges to $x \in X$ if
  and only if $\Diamond x = \alpha_{\mathcal{F}}$.
\end{lem}
\proof
`If' part is obvious from the definition of
$\alpha_{\mathcal{F}}$.

For `only if' part, suppose that 
$\mathcal{F}$ converges to $x \in X$. Let $a = \ball_{d}(y,\delta) \in
\Diamond x$. Choose $\theta \in \PRat$ such that $d(y,x) + \theta <
\delta$. Then, $\ball_{d}(x,\theta) <_X a$. Since
$\mathcal{F}$ converges to $x$, we have $\Ball_{d}(x,\theta) \in
\mathcal{F}$, and so $a \in \alpha_{\mathcal{F}}$. Conversely, 
let $a = \ball_{d}(y,\delta) \in \alpha_{\mathcal{F}}$. Then, there
exists $\ball_{d'}(y',\delta') <_{X} a$ such that 
$\Ball_{d'}(y',\delta') \in \mathcal{F}$. Choose $\theta \in \PRat$
such that $d(y,y') + \delta' + \theta < \delta$. Since
$\mathcal{F}$ converges to $x$, we have $\Ball_{d}(x,\theta) \in
\mathcal{F}$. Since $\mathcal{F}$ is a filter, we have 
$d(y,x) \leq d(y,y') + d(y',x) <d(y,y')  + \delta' + \theta < \delta$.
Hence $a \in \Diamond x$.
\qed

\begin{prop}\label{prop:CompChar}
  A symmetric gus $(X,M)$ is complete if and only if the function
  $i_{X} \colon X \to \Pt(\mathcal{U}(X))$ is a surjection, i.e.
  for each $\alpha \in \Pt(\mathcal{U}(X))$ there exists $x \in X$ such
  that $\alpha = \Diamond x$.
\end{prop}
\proof
By Lemma \ref{lem:FiltPt} and Lemma \ref{lem:FiltConv}.
\qed

\begin{lem}[cf.\ {\cite[Theorem
  2.7]{PalmgrenLocalicCompletion}}]\label{lem:DenseIso}
  Let $(X,M)$ and $(Y,N)$ be symmetric gus's where $M = \left\{
    d_{i} \mid i \in I \right\}$ and $N = \left\{ \rho_{i} \mid i \in I
  \right\}$ are indexed by the same set $I$, and which satisfies
  $d_{i} \leq d_{j} \iff \rho_{i} \leq \rho_{j}$ for each $i,j \in I$.
   Let $f \colon X \to Y$ be an isometry in the sense that 
   $d_{i} = \rho_{i} \circ f \times f$ for each $i \in I$
   with a dense image. Then, the functor $\overline{(-)} \colon \GUS
   \to \FTop$ sends $f$ to an isomorphism $\overline{f} \colon
   \mathcal{U}(X) \to \mathcal{U}(Y)$.
\end{lem}
\proof
The reader is referred to a  quite similar proof of \cite[Theorem
2.7]{PalmgrenLocalicCompletion}. 
Note that 
\begin{equation}
  \ball_{d_{i}}(x,\varepsilon) \mathrel{\overline{f}}
  \ball_{\rho_{j}}(y,\delta)
  \iff
  \ball_{\rho_{i}}(f(x),\varepsilon) <_{Y} 
  \ball_{\rho_{j}}(y,\delta).
  \tag*{$\qEd$}
\end{equation}

By Proposition \ref{prop:DenseEmb} and Lemma \ref{lem:DenseIso}, the isometry $i_{X} \colon X \to
\Pt{\mathcal{U}(X)}$ gives rise to an isomorphism
$\overline{i_{X}} \colon \mathcal{U}(X) \to
\mathcal{U}(\Pt(\mathcal{U}(X)))$ defined by 
\[
  \ball_{d}(x,\varepsilon) \mathrel{\overline{i_{X}}}
  \ball_{\widetilde{\rho}}(\alpha,\delta) \iff
  \ball_{\widetilde{d}}(\Diamond x, \varepsilon)
  <_{\Pt(\mathcal{U}(X))} \ball_{\widetilde{\rho}}(\alpha,\delta).
\]
Applying the operation $\Pt(-)$ to $\overline{i_X}$, we obtain an
isomorphism $\Pt(\overline{i_{X}}) \colon\Pt(\mathcal{U}(X)) \to
\Pt(\mathcal{U}(\Pt(\mathcal{U}(X))))$.

\begin{lem}
  For each $\alpha \in \Pt(\mathcal{U}(X))$, we have
  \[
    \Pt(\overline{i_{X}})(\alpha) = \Diamond \alpha.
  \]
\end{lem}
\proof
First, suppose $\ball_{\widetilde{d}}(\beta, \delta) \in
\Pt(\overline{i_{X}})(\alpha)$. Then, there exists
$\ball_{\rho}(x,\varepsilon) \in \alpha$ such that
$\ball_{\widetilde{\rho}}(\Diamond x, \varepsilon)
<_{\Pt(\mathcal{U}(X))} \ball_{\widetilde{d}}(\beta, \delta)$.
Then, 
  $
  \widetilde{d}(\beta, \alpha)  \leq
  \widetilde{d}(\beta,\Diamond x) + \widetilde{d}(\Diamond x, \alpha)
  < \widetilde{d}(\beta, \Diamond x) + \varepsilon < \delta.
  $
Hence, $\ball_{\widetilde{d}}(\beta, \delta) \in \Diamond \alpha$.

Conversely, let $\ball_{\widetilde{d}}(\beta, \delta) \in \Diamond \alpha$.
Choose $\theta \in \PRat$ such that  $\widetilde{d}( \beta, \alpha) + 2
\theta < \delta$.
Since $i_{X} \colon X \to \Pt(\mathcal{U}(X))$ is dense, there exists $x
\in X$ such that $\widetilde{d}(\alpha, \Diamond x) < \theta$ so that
$\ball_{d}(x,\theta) \in \alpha$. Moreover, 
  $
  \widetilde{d}(\beta, \Diamond x) + \theta \leq 
  \widetilde{d}(\beta, \alpha) + \widetilde{d}(\alpha, \Diamond x) +
  \theta < 
  \widetilde{d}(\beta, \alpha) + 2\theta < \delta.
  $
  Thus, $\ball_{\widetilde{d}}(\Diamond x,\theta)
<_{\Pt(\mathcal{U}(X))} \ball_{\widetilde{d}}(\beta,\delta)$, and
hence $\ball_{\widetilde{d}}(\beta, \delta) \in
\Pt(\overline{i_{X}})(\alpha)$.
\qed

By Proposition \ref{prop:SeparatedGUS} and Proposition
\ref{prop:CompChar}, $\left(\Pt(\mathcal{U}(X)),
\widetilde{M}\right)$ is a complete separated symmetric gus. Thus, we
conclude as follows.
\begin{thm}\label{thm:Completion}
  The isometry $i_{X} \colon  (X,M) \to \left(\Pt(\mathcal{U}(X)),
  \widetilde{M}\right)$ is a completion of $(X,M)$.
\end{thm}

\begin{rem}\label{rem:Predicativity}
  By Theorem \ref{thm:ClosedEmbedding},
  the symmetric gus
  $\Pt(\mathcal{U}(X,M))$ embeds into the product
  $\Pt(\prod_{d \in M}\mathcal{M}(X,d)) \cong \prod_{d \in
  M}\Pt(\mathcal{M}(X,d))$ of complete generalised uniform spaces as a
  closed subclass. This is one
  of the well-known construction of completions of uniform spaces
  \cite{KelleyGenTop}. 
  In particular, for a finite Dedekind symmetric gus $(X,M)$,
  the class $\Pt(\mathcal{U}(X,M))$ gives the usual completion of $(X,M)$
  with Cauchy filters.
  
  The construction of $\Pt(\mathcal{U}(X,M))$, however, is problematic from a
  predicative point of view because we have yet to prove that
  $\Pt(\mathcal{U}(X,M))$ is a set.
  If $(X,M)$ is a finite Dedekind symmetric gus, this problem can be
  addressed by using the fact that
  for each pseudometric $d \in M$, the class $\Pt(\mathcal{M}(X,d))$
  is isometric to the usual metric completion of $(X,d)$ with Cauchy
  sequences \cite[Theorem 2.4]{PalmgrenLocalicCompletion} (this
  requires $\ACw$). The completion of $(X,M)$ is then obtained as a
  closed subset of the product $\prod_{d \in M}\Pt(\mathcal{M}(X,d))$.
  Hence, the construction of the completion of Bishop's notion of
  uniform space \cite[Chapter 4, Problem 17]{Bishop-67} is predicative
  and unproblematic under the assumption of $\ACw$, which is a usual
  practice in Bishop constructive mathematics. See below for another way to
  cope with this predicativity issue.
\end{rem}
\subsubsection{Predicativity issues}\label{sec:Predicative}
We present another way to cope with the predicative issue
raised in Remark \ref{rem:Predicativity}. 
This method is specific to $\CZF$, but it avoids $\ACw$ and
works for any symmetric gus.
Note that to show that $\Pt(\mathcal{U}(X,M))$ forms a set,
it suffices to show that for each symmetric gms $(X,d)$, the class
$\Pt(\mathcal{M}(X,d))$ forms a set.

Recall that Fullness is the following statement, which is valid in $\CZF$.
  \[
     \forall X  \forall Y 
     \exists R 
    \left[ R \subseteq \mv(X,Y)  \amp  \forall s \in
      \mv(X,Y) \exists\, r \in R \; r \subseteq
    s\right]
  \]
where $\mv(X,Y)$ is the \emph{class} of total relations between
sets $X$ and $Y$.

Now, let $(X,d)$ be a symmetric gms. By Fullness, there exists a
subset 
\[
  R \subseteq \mv(\PRat, U_{(X,d)})
\]
such that for each $s \in \mv(\PRat, U_{(X,d)})$ there exists $r \in
R$ such that $r \subseteq s$.
Define a set
\[
  \overline{R}
  \defeql
  \left\{ r \in R \mid \left( \forall (\varepsilon, a)  \in r \right)
    a \in \mathcal{C}^{\varepsilon}_{d} \amp \left(
    \forall (\varepsilon,a),(\delta,b) \in r \right) a_{*} \meets
    b_{*} \right\}.
\]
The following is obvious.
\begin{lem}\label{lem:alphar}
  For each $r \in \overline{R}$, the subset
  \[
    \alpha_{r} \defeql \left\{ a \in U_{(X,d)} \mid \left( \exists
    (\varepsilon, b) \in r \right) b <_{X} a\right\}
  \]
  is in $\Pt(\mathcal{M}(X,d))$.
\end{lem}
Conversely, given $\alpha \in \Pt(\mathcal{M}(X,d))$, put
\[
  r_{\alpha} \defeql \left\{ (\varepsilon, a) \in \PRat \times
    U_{(X,d)} \mid a \in \alpha \cap
    \mathcal{C}^{\varepsilon}_{d} \right\}.
\]
By \ref{U2}, we have $r_{\alpha} \in \mv(\PRat, U_{(X,d)})$. Hence, by
Fullness there exists $r \in R$ such that $r \subseteq r_{\alpha}$.
Clearly $r \in \overline{R}$, and $\alpha_{r} \subseteq \alpha$.
Let $a \in \alpha$, and write $a = \ball_{d}(x,\varepsilon)$. By
\ref{U1}, there exists $\delta < \varepsilon$ such that
$\ball_{d}(x,\delta) \in \alpha$. Choose $\theta \in \PRat$ such that
$\delta + 2 \theta < \varepsilon$. By \ref{U2}, there exists
$\ball_{d}(y,\theta) \in \alpha_{r}$, and hence by \ref{P2} applied to
$\alpha$, we have $d(x,y) < \delta + \theta$. Then, $d(x,y) + \theta <
\delta + 2 \theta < \varepsilon$, and so $\ball_{d}(y,\theta) <_{X}
\ball_{d}(x,\varepsilon) = a$. Then, $a \in \alpha_{r}$ by \ref{P3}.
Thus $\alpha_{r} = \alpha$.

\begin{prop}
  The assignment $r \mapsto \alpha_{r}$ is a surjection from
  $\overline{R}$ to $\Pt(\mathcal{M}(X,d))$.
  Hence, $\Pt(\mathcal{M}(X,d))$ is a set.
\end{prop}
The argument following Lemma \ref{lem:alphar} shows that every formal
point of $\Pt(\mathcal{M}(X,d))$ is maximal. Hence, the above
predicativity issue is a special case of the example discussed in
\cite{PalmgrenMaxPartialPt}, \cite{vandenBergNID}, and
\cite{SGAxiom}. Note that we do not require any extra axiom of $\CZF$
discussed in \cite{vandenBergNID} and \cite{SGAxiom}.
  
\section{Functorial embedding of locally compact uniform spaces}\label{sec:FunctEmb}
In this section, we work with separated finite Dedekind symmetric
gus's, which we simply call uniform spaces. This is the notion of
uniform space treated in \cite{Bishop-67}.
We show that the category of locally compact uniform spaces can be
embedded into that of overt locally compact completely regular formal topologies
by extending the construction of localic completions to a full and
faithful functor. 
The results in this section are rather straightforward generalisations
of the corresponding results for metric spaces  by Palmgren
\cite[Section 4]{PalmgrenLocalicCompletion}. However, our
characterisation of the cover of the localic completion of a locally
compact uniform space (and thus also that of a locally compact metric space) may be interesting.
We believe that our characterisation is an improvement over the
previous one for metric spaces \cite[Theorem 4.17]{PalmgrenLocalicCompletion}.

\subsection{Complete regularity}\label{sec:RegularCompact}
Classically, the topology associated with a uniform space is
completely regular. Hence, the localic completion of a uniform space
should be (point-free) completely regular. Before making it precise,
we recall the relevant notions from \cite{Curi200385}.

Let $\mathcal{S}$ be a formal topology. For a subset
$U \subseteq S$, let
  $
  U^{*} \defeql \left\{ b \in S \mid b \downarrow U \cov
  \emptyset \right\},
  $
and for any two subsets $U,V \subseteq S$, write $U \lll V$ if $S
\cov U^{*} \cup V$. We write $a \lll b$ for $\left\{ a
\right\} \lll \left\{ b \right\}$.
Put $\URat \defeql \left\{ q \in \Rat \mid 0 \leq q \leq 1 \right\}$.
For subsets $U,V \subsets S$, a \emph{scale} from $U$ to $V$ is a family
$\left ( U_q\right)_{q \in \URat}$ of subsets of $S$ such that $U \cov
U_0$, $U_1 \cov V$ and $U_p \lll U_q$ for each $p,q \in \URat$ such
that $p < q$.
We write $U \creg V$ if there exists a scale
from $U$ to $V$ and write $a \creg b$ for $\left\{
 a \right\} \creg \left\{ b \right\}$.

\begin{defi}\label{def:CompEnumCReg}
  A formal topology $\mathcal{S}$ is \emph{completely
  regular} if it is equipped with a function $\rc \colon S \to \Pow(S)$ such
that
  \begin{enumerate}
    \item\label{eq:rc1} $\bigl( \forall b \in \rc(a) \bigr)\: b \creg a$, 
    \item\label{eq:rc2} $a \cov \rc(a)$ 
  \end{enumerate}
   for all $a \in S$.
\end{defi}

\begin{lem}\label{lem:LCompReg}
  For any symmetric Dedekind gus $(X,M)$, we have
  \[
    a <_X b \implies a \lll b
  \]
  for all $a,b \in U_X$.
\end{lem}
\proof
  Let $a,b \in U_X$, and suppose that $a <_X b$. Write
  $a = \ball_{d}(x,\varepsilon)$, and choose $\theta \in \PRat$
  such that $\ball_{d}(x,\varepsilon + 3\theta) <_X b$.
  Let $c = \ball_{d}(z,\theta) \in \mathcal{C}^{\theta}_{d}$.
  Then, either $d(x,z) > \varepsilon + \theta$ or 
  $d(x,z) < \varepsilon + 2\theta$.
  In the former case, for any $c' \in a \downarrow c$,
  we have $d(x,z) < \varepsilon + \theta$, a contradiction.
  Thus, $a \downarrow c \cov_{X} \emptyset$ and so  $c \in a^{*}$. In the latter case, we have
  \[
    d(x,z) + \theta \leq \varepsilon + 3\theta,
  \]
  so $\ball_{d}(z,\theta) \leq_X \ball_{d}(x,\varepsilon + 3
  \theta) <_X b$. Hence, $U_X \cov_{X} a^{*} \cup \left\{ b
  \right\}$ by \ref{U2}. Therefore, $a \lll b$.
\qed

\begin{prop}
  The localic completion of a symmetric Dedekind gus is completely regular.
\end{prop}
\proof
Let $(X,M)$ be a symmetric Dedekind gus. For each
$\ball_{d}(x,\varepsilon) \in U_{X}$, put
\[
  \rc(\ball_{d}(x,\varepsilon)) \defeql \left\{
    \ball_{d}(x, \delta) \in U_X \mid \delta < \varepsilon \right\}.
\]
By \ref{U1}, we have $a \cov_{X} \rc(a)$ for each $a \in U_{X}$.
Let $a,b \in U_X$ such that $b \in \rc(a)$, and write $a =
\ball_{d}(x,\varepsilon)$ and $b = \ball_{d}(x,\delta)$. 
Since $\delta < \varepsilon$, there exists an order preserving bijection
between $[\delta,\varepsilon] \cap \Rat$ and $\URat$. Thus, $b \creg a$
by Lemma \ref{lem:LCompReg}, and so $\mathcal{U}(X)$ is completely
regular.
\qed
%

\subsection{Compactness and local compactness}\label{sec:CompLComp}
The following notion of local compactness generalises the
corresponding notion for metric spaces by Bishop \cite[Chapter 4, Definition 18]{Bishop-67}. 
\begin{defi}\label{def:LKUSpa}
  A uniform space $(X,M)$ is \emph{totally bounded} if for each
  $d \in M$, the pseudometric space $(X,d)$ 
  is totally bounded, i.e.\
  \begin{align*}
    \left( \forall \varepsilon \in \PRat \right)
    \left( \exists Y_{\varepsilon} \in \Fin(X) \right)
    X \subseteq \bigcup_{y \in Y_{\varepsilon}} \Ball_{d}(y,\varepsilon).
  \end{align*}
  The set $Y_{\varepsilon}$ is called an \emph{$\varepsilon$-net} to
  $X$ with respect to $d$. 
  A uniform space is \emph{compact} if it is complete and totally bounded.
  A uniform space $X$ is \emph{locally compact} if for each open
  ball $\Ball_{d}(x,\varepsilon)$ of $X$, there exists a compact subset $K \subseteq X$
  such that $\Ball_{d}(x,\varepsilon) \subseteq K$.
  Thus, every compact uniform space is locally compact.
\end{defi}
\begin{rem}
The above notion of local compactness is not invariant under
homeomorphisms, which is often considered unsatisfactory. For
example, the real line is locally compact in the above sense, but
$(0,1)$ is not. In this paper, we are not aiming at a better
alternative to Bishop's definition of local compactness but one that
is compatible with it.
\end{rem}

Compactness in formal topology is defined by the covering compactness.
\begin{defi}
A formal topology $\mathcal{S}$ is \emph{compact} if 
\[
S \cov U \implies \bigl( \exists U_0  \in \Fin(U) \bigr)\: S\cov U_0
\]
for all $U \subseteq S$.
\end{defi}

\begin{defi}
  Let $\mathcal{S}$ be a formal topology. For each $a, b \in S$ define
  \begin{equation*}\label{form:waybelow}
     a \ll b \defeqiv \bigl( \forall U \in \Pow(S) \bigr)\left[ \:
     b \cov U
     \implies \bigl( \exists U_0 \in \Fin(U) \bigr)\; a \cov U_0
   \right].
  \end{equation*}
  A formal topology $\mathcal{S}$ is \emph{locally compact} if it is
  equipped with a function $\wb \colon S \to \Pow(S)$ such that
  \begin{enumerate}
    \item\label{eq:wb1} $\bigl( \forall b \in \wb(a) \bigr)\: b \ll a$, 
    \item\label{eq:wb2} $a \cov \wb(a)$
  \end{enumerate}
 for all $a \in S$.
Since the relation $\ll$ is a proper class in general, the function
$\wb \colon S \to \Pow(S)$ is an essential part of the definition of
local compactness.
\end{defi}

Given a uniform space $(X,M)$, define a relation $\sqsubseteq$ on
$\Pow(U_X)$ by
\begin{align*}\label{def:Ballwise}
  U \sqsubseteq V \defeqiv \left( \exists (d,\varepsilon) \in \Ent(M) \right)
  U \downarrow \mathcal{C}_{d}^{\varepsilon} \leq_{X} V
\end{align*}
for all $U,V \subseteq U_X$. By \ref{U2}, we have
\begin{align*}
  U \sqsubseteq V \implies U \cov_{X} V.
\end{align*}

\begin{lem}\label{lem:Sqeq}
  Let $(X,M)$ be a uniform space. Then,
  the following are equivalent for all $a \in U_X$ and $U \subseteq U_X$:
  \begin{enumerate}
    \item\label{lem:Sqeq1} $\left( \exists V \in \Fin(U_X) \right) a_{*} \subseteq
      V_{*} \amp V <_X U$;
    \item\label{lem:Sqeq2} $\left( \exists V \in \Fin(U_X) \right) a \sqsubseteq V  <_X  U$;
    \item\label{lem:Sqeq3} $\left( \exists V \in \Fin(U_X) \right) a \cov_X V  <_X  U$.
  \end{enumerate}
\end{lem}
\proof
  \eqref{lem:Sqeq1} $\Rightarrow$ \eqref{lem:Sqeq2} Suppose that
  \eqref{lem:Sqeq1} holds. Then, there exists
  $V \in \Fin(U_X)$ such that $a_{*} \subseteq V_{*}$  and $V <_X U$.
  Write $V = \left\{ \ball_{d_0} (x_0,
  \varepsilon_0),\dotsc,\ball_{d_{n}} (x_{n}, \varepsilon_{n})
  \right\}$, and choose $\theta \in \PRat$ such that $V' \defeql \left\{
  \ball_{d_{i}}(x_i,\varepsilon_i + \theta) \mid i \leq n\right\} <_X U$.
  Let $d \defeql \sup\left\{d_{i} \mid i \leq n\right\}$, and let
  $\ball_{\rho}(y,\delta) \in a
  \downarrow \mathcal{C}^{\theta}_{d}$. Then, there exists $i \leq n$
  such that $d_i(x_i,y) < \varepsilon_i$, so we have 
    $
    d_i(x_i,y) + \delta < \varepsilon_{i} + \delta \leq
    \varepsilon_{i} + \theta.
    $
  Thus, $\ball_{\rho}(y,\delta) <_X \ball_{d_{i}}(x_i, \varepsilon_i +
  \theta)$, and hence $a \sqsubseteq V'$.

  \eqref{lem:Sqeq2} $\Rightarrow$ \eqref{lem:Sqeq3} We have
  $a \sqsubseteq V \implies a \cov_{X} V$.

  \eqref{lem:Sqeq3} $\Rightarrow$ \eqref{lem:Sqeq1} We have
  $a \cov_{X} V \implies a_{*} \subseteq V_{*}$.
\qed

The  cover of the localic completion of a locally compact uniform space
admits an elementary characterisation.
\begin{prop}\label{lem:CoverLCM}
  Let $(X,M)$ be a locally compact uniform space. Then, the
  following are equivalent for all $a \in U_X$ and $U \subseteq U_X$:
  \begin{enumerate}
    \item\label{lem:CoverLCM1} $a \cov_X U$;
    \item\label{lem:CoverLCM2} $\left(\forall b <_X a \right)\left( \exists V \in \Fin(U_X)\right) b_{*} \subseteq V_{*} \amp  V <_X U$;
    \item $\left(\forall b <_X a \right)\left( \exists V \in \Fin(U_X)\right) b \sqsubseteq V <_X U$;
    \item $\left(\forall b <_X a \right)\left( \exists V \in \Fin(U_X)\right)
      b \cov_{X} V <_X U$.
  \end{enumerate}
\end{prop}
\proof
  By Lemma \ref{lem:Sqeq}, it suffices to show that \eqref{lem:CoverLCM1} implies \eqref{lem:CoverLCM2}.
  Given $U \subseteq U_X$, define a predicate $\Phi_{U}$ on $U_X$ by
  \begin{align*}
    \Phi_{U}(a) \defeqiv \left( \forall b <_X a \right) \left( \exists V
    \in \Fin(U_X)\right) b_{*} \subseteq V_{*} \amp V <_X U.
  \end{align*}
  We show that 
  \[
  a \cov_X U \implies \Phi_{U}(a)
  \]
  for all $a \in U_X$ by induction on $\cov_{X}$.
  We must check the conditions \ref{ID1} -- \ref{ID3}
  for the localised axiom-set consisting of \ref{U1} and \ref{U2'}.

  The conditions \ref{ID1} and \ref{ID2} are straightforward
  to check, using Lemma \ref{lem:misc}. For \ref{ID3},
  we have two axioms to be checked.

  \ref{U1} $\displaystyle\frac{\left( \forall b <_X a \right)
  \Phi_{U}(b)}{\Phi_{U}(a)}$: Suppose that $\Phi_{U}(b)$ for all
  $b <_X a$. Let $b <_X a$. Then, there exists $c \in U_X$ such that
  $b <_X c <_X a$. Since $\Phi_{U}(c)$, there exists $V \in \Fin(U_X)$
  such that $b_{*} \subseteq V_{*}$ and $V <_X U$.
  Hence, $\Phi_{U}(a)$.

  \ref{U2'} $\displaystyle
  \frac{\left( \forall c \in  \mathcal{C}^{\theta}_{\rho} \downarrow a \right)
  \Phi_{U}(c)}{\Phi_{U}(a)}$ for
  each $(\rho,\theta) \in \Ent(M)$:
  Suppose that $\Phi_{U}(c)$ for all $c \in \mathcal{C}^{\theta}_{\rho}
  \downarrow a $.
  Let $b <_X a$, and write $a = \ball_{d}(x,\varepsilon)$ and $b =
  \ball_{d'}(y,\delta)$.
  Since $(X,M)$ is locally compact, there exists a compact
  subset $K \subseteq X$ such that $\Ball_{d}(x,\varepsilon) \subseteq K$.
  Choose $\xi \in \PRat$ such that $2 \xi < \theta$ and $d(x,y) + \delta + 4\xi <
  \varepsilon$. Let $Z \defeql \left\{ z_{0},\dotsc, z_{n-1} \right\}$ be a
  $\xi$-net to $K$ with respect to $\rho' \defeql \sup\left\{ d, \rho \right\}$.
  Split $Z$ into two finitely enumerable subsets $Z^{+}$ and $Z^{-}$ such
  that $Z = Z^{+} \cup Z^{-}$ and
  \begin{itemize}
    \item  $z \in Z^{+}  \implies d(z,x) < \varepsilon - 2\xi$,
    \item  $z \in Z^{-} \implies d(z,x) > \varepsilon - 3\xi$.
  \end{itemize}
  Let $z \in Z^{+}$. Since $\ball_{\rho'}(z,2\xi)
  \in\mathcal{C}^{\theta}_{\rho}  \downarrow  a$,
  we have $\Phi_{U}(\ball_{\rho'} (z,2\xi))$.  Hence
  there exists $V_z \in \Fin(U_X)$ such that $\Ball_{\rho'}(z,\xi) \subseteq
  {V_{z}}_{*}$ and $V_z <_X U$. Since $Z^{+}$ is finitely enumerable,
  there exists $V \in \Fin(U_X)$ such that 
    $ \bigcup_{z \in Z^{+}} \Ball_{\rho'}(z, \xi) \subseteq V_{*} $
  and $V <_X U$.
  Now, it suffices to show that $b_{*} \subseteq \bigcup_{z \in
  Z^{+}}\Ball_{\rho'}(z, \xi)$. Let $y' \in b_{*}$. Then, 
  there exists $i < n$ such that $\rho'(y',z_{i}) < \xi$.
  Then
  \begin{align*}
    d(z_i,x) 
    &\leq d(z_{i},y') + d(y',y) + d(y,x) \\
    &< \xi + \delta + d(y,x) < \varepsilon - 3 \xi,
  \end{align*}
  and  thus $z_{i} \in Z^{+}$. Hence, $y' \in \bigcup_{z \in
  Z^{+}}\Ball_{\rho'}(z, \xi)$, and therefore $\Phi_{U}(a)$.
\qed 

\begin{rem}\label{rem:REA}
  By Proposition \ref{lem:CoverLCM}, inductive generation of the cover
  of the localic completion of a locally compact uniform space
  does not require the Regular Extension Axiom.
\end{rem} 

\begin{example}[{\cite[Example 3.3]{PalmgrenLocalicCompletion}}]
\label{eg:R2}
  Consider the real plane $X \defeql \Real^{2}$, which is a locally compact
  uniform space (it is even a metric space). Let $x = (0,0),
  y = (-4,0)$, and $z = (4,0)$. Put $a = \ball_{d}(x,3)$, $b =
  \ball_{d}(y,5)$, and $c = \ball_{d}(z,5)$, where $d$ is the standard
  distance on $\Real^{2}$ (see the left figure in Figure
  \ref{fig:egR2}). Then, $a \cov_{X} \left\{b,c\right\}$.
  This can be seen as follows: if we shrink the radius of $a$ by
  $\varepsilon \in \PRat$ and let $a' = \ball_{d}(x,3 - \varepsilon)$,
  we can find a sufficiently small $\delta \in \PRat$ such that $a'
  \downarrow \mathcal{C}_{d}^{\delta} \leq_{X} \left\{b,c \right\}$ as
  can be visually seen from the right figure in Figure \ref{fig:egR2}.
  \begin{figure}[h]
    \begin{tikzpicture}[scale=1.2,transform shape, every
      node/.style={font=\scriptsize}]
      \draw (-1, 0) circle (1.25);
      \draw ( 1, 0) circle (1.25);
      \draw ( 0, 0) circle (0.75);
      \node at (0,0.32) [inner sep=0pt,circle,fill=white] {$a$};
      \node at (-1.2,0.5) {$b$};
      \node at (1.2,0.5) {$c$};
      \fill (-1, 0) circle (1.0pt) ;
      \fill (1, 0) circle (1.0pt);
      \fill (0, 0) circle (1.0pt);
      \draw[{Latex}-{Latex}] (0,0) -- (-0.60,-0.45)
      node[circle, inner sep=0pt, midway, fill=white] {$3$};
      \draw[{Latex}-{Latex}] (-1,0) -- (-2,-0.75)
      node[midway,circle,inner sep=0pt, fill=white] {$5$};
      \draw[{Latex}-{Latex}] (1,0) -- (2,-0.75)
      node[midway,circle,inner sep=.1pt, fill=white] {$5$};
      \draw[{Latex}-{Latex}] (0,0) -- (1,0) node[midway,circle,inner sep=0pt, fill=white] {$4$};
    \end{tikzpicture}
    \qquad
    \qquad
    \begin{tikzpicture}[scale=1.2,transform shape, every
      node/.style={font=\scriptsize}]
      \draw (-1, 0) circle (1.25);
      \draw ( 1, 0) circle (1.25);
      \draw ( 0, 0) circle (0.60);
      \node at (0.02,0.3) [inner sep=0pt,circle,fill=white] {$a'$};
      \node at (-1.2,0.5) {$b$};
      \node at (1.2,0.5) {$c$};
      \fill (-1, 0) circle (1.0pt) ;
      \fill (1, 0) circle (1.0pt);
      \fill (0, 0) circle (1.0pt);
      \draw[{Latex}-{Latex}] (0,0) -- (-0.45,-0.40);
      \node[inner sep=0pt] (A) at (1.0,-0.8) {$3 - \varepsilon$};
      \path[draw] (-.225,-0.2) .. controls (0.1,-.50) and  (0.5,-.03) .. (A);
    \end{tikzpicture}
    \caption{Example \ref{eg:R2}}
    \label{fig:egR2}
  \end{figure}
  Note that we should not conclude $a \cov_{X} \left\{b,c\right\}$ from
  the left figure just because $a_{*} \subseteq b_{*} \cup c_{*}$.
  This relies on the spatiality of $\mathcal{U}(\Real^{2})$.
\end{example}

\begin{example}\label{eg:FR}
  The characterisation of the cover in Proposition \ref{lem:CoverLCM}
  can be considered as a natural generalisation of the characterisation
  of the cover of the localic reals by Johnstone \cite[Chapter IV,
  Section 1.1, Lemma]{johnstone-82}. We restate his characterisation
  in terms of formal reals $\FR$, which is shown to be identical
  to the localic completion of the space of rational numbers 
  \cite[Example 2.2]{PalmgrenLocalicCompletion}.
  Recall that $\FR$ is a formal topology with a base $S_{\FR} = \left\{
  (p,q) \in \Rat \times \Rat \mid p < q \right\}$ ordered by $(p,q)
  \leq_{\FR}
  (r,s) \defeqiv r \leq p \amp q \leq s$.  For each $(p,q), (r,s) \in
  S_{\FR}$ define
    $
    (p,q) <_{\FR} (r,s)\defeqiv r < p < q < s.
    $
  The axioms of $\FR$ are the following:
    \begin{enumerate}[label=(R\arabic*)]
    \item\label{R1} $(p,q) \cov_{\FR} \left\{ (r,s) \in S_{\FR} \mid
      (r,s) <_{\FR} (p,q)    \right\}$,
    \item\label{R2} $(p,q) \cov_{\FR} \left\{ (p,s), (r,q) \right\}$ \
      for each $p < r < s < q$.
  \end{enumerate}
 Then, for each $U \subseteq S_{\FR}$ and $(p,q) \in S_{\FR}$, we have
    \begin{align*}
      (p,q) \cov_{\FR} U \iff& 
      \forall (p',q') <_{\FR} (p,q) \, \exists (p_{0}, q_{0}), \dots,
      (p_{n}, q_{n})  \in {\downarrow} U  \\
      & \qquad p' = p_{0} \amp q_{n} = q'  \amp
      \left( \forall i < n \right) p_{i} \leq p_{i+1} < q_{i} \leq
      q_{i+1},
    \end{align*}
  where ${\downarrow} U$ is the downward closure of $U$ with respect
  to $\leq_{\FR}$. The proof is by a straightforward induction on
  $\cov_{\FR}$.
\end{example}

\begin{cor}\label{col:StrictOrderWayBelow}
  For any locally compact uniform space $X$, we have
  \begin{align*}
    a <_X b \implies a \ll b
  \end{align*}
  for all $a,b \in U_X$.
\end{cor}

By \ref{U1}, we obtain the following.
\begin{thm}\label{thm:LCompLCMLK}
 The localic completion of a locally compact uniform space is locally compact.
\end{thm}

\begin{thm}\label{thm:LCompTBK}
  A uniform space $(X,M)$ is totally bounded if and only if $\mathcal{U}(X)$ is compact.
\end{thm}
\proof
  Suppose that $(X,M)$ is totally bounded. Let $U \subseteq U_X$, and suppose
  that $U_X \cov_X U$. Choose any $d \in M$ and $\varepsilon \in \PRat$,
  and let $\left\{ x_{0},\dotsc, x_{n-1} \right\}$ be an
  $\varepsilon$-net to $X$ with respect to $d$.
  By \ref{U2}, we have $U_X \cov_X \mathcal{C}^{\varepsilon}_{d}
  \cov_X \left\{ \ball_{d}(x_{i},2\varepsilon) \mid i < n\right\}$.
  Thus, there exists $\ball_{d}(y,\delta) \in U_X$ such that $U_X \cov_X
  \ball_{d}(y,\delta)$. Since $\ball_{d}(y,2\delta) \cov_X U$, 
  there exists $V \in \Fin(U_X)$ such that $\ball_{d}(y,\delta) \cov_X V
  <_X U$ by Proposition \ref{lem:CoverLCM}. Then, there exists $U_{0}\in
  \Fin(U)$ such that $U_X \cov_X U_0$.  Therefore, $\mathcal{U}(X)$
  is compact.

  Conversely, suppose that $\mathcal{U}(X)$ is compact. Let $d
  \in M$ and $ \varepsilon \in \PRat$. Since $U_X \cov_X
  \mathcal{C}^{\varepsilon}_{d}$, there exists
  $V = \left\{ \ball_{d}(x_0,\varepsilon),\dotsc,
       \ball_{d}(x_{n-1},\varepsilon) \right\}
           \in \Fin(\mathcal{C}^{\varepsilon}_{d})$
  such that $U_X \cov_X V$. Then, $X = {U_X}_{*} \subseteq
  V_{*} = \bigcup_{i < n} \Ball_{d}(x_i,\varepsilon)$, and hence
  $\left\{ x_0,\dotsc, x_{n-1} \right\}$ is an $\varepsilon$-net to
  $X$ with respect to $d$. Therefore, $X$ is totally bounded.
\qed

\subsection{Functorial embedding}\label{subsec:FunctEmb}
\begin{defi}\label{def:ContFun}
%
  A function $f \colon X \to Y$ from a  locally compact uniform space
  $(X,M)$ to a uniform space  $(Y,N)$ is \emph{continuous} if $f$ is
  uniformly continuous on each open ball of $X$, i.e.
  for each $x \in X$, $d \in M$ and $\varepsilon \in \PRat$, 
  \begin{multline*}
    \left( \forall \rho \in N \right) \left( \forall \delta \in
    \PRat \right) \left( \exists d' \in M \right) \left( \exists
    \varepsilon' \in \PRat \right) \\
    \left[ \left( \forall x_{1},x_{2} \in \Ball_{d}(x,\varepsilon) \right)
    d'(x_{1},x_{2}) < \varepsilon' \implies
  \rho(f(x_{1}),f(x_{2})) < \delta‘\right].
  \end{multline*}
\end{defi}


Since the image of a totally bounded uniform
space under a uniformly continuous function is again totally bounded,
continuous functions between locally compact uniform spaces are closed
under composition.  Thus, the locally compact uniform
spaces and continuous functions form a category,
which we denote by $\LKUSpa$.

\begin{lem}
  A locally compact uniform space is complete.
\end{lem}
\proof
  Let $(X,M)$ be a locally compact uniform space.  Let
  $\mathcal{F}$ be a Cauchy filter on $X$.  Choose any  $d \in M$ and
  $\varepsilon \in \PRat$.  By \ref{CF4}, there exists $x \in X$
  such that $\Ball_{d}(x,\varepsilon) \in \mathcal{F}$.
  Since $X$ is locally compact,
  there exists a compact subset $K \subseteq X$ such that
  $\Ball_{d}(x,\varepsilon) \subseteq K$.
  Let $\mathcal{G} = \left\{ U \in \mathcal{F}
  \mid U \subseteq K\right\}$. It is easy to see that
  $\mathcal{G}$ is a Cauchy filter on $K$. 
  Since $K$ is complete, $\mathcal{G}$ converges to some $z
  \in K$. Since $\mathcal{G} \subseteq \mathcal{F}$, $\mathcal{F}$ also
  converges to $z$. 
\qed
Thus, for each locally compact uniform space $X$, the embedding $i_{X}
: X \to \Pt(\mathcal{U}(X))$ defined by \eqref{eq:Completion} is a uniform isomorphism.

Given any function $f \colon X \to Y$ between uniform spaces $X$ and
$Y$, define a relation  $r_{f} \subseteq U_X \times U_Y$ by
\begin{equation} \label{def:UMap}
  a \mathrel{r_{f}} b \defeqiv \left( \exists b' <_Y b \right)
  f[a_{*}]\subseteq b'_{*}
\end{equation}
for each $a \in U_X$ and $b \in U_{Y}$.
The relation $r_{f}$ is a natural generalisation of the
relation $D_{f}$ defined by Palmgren in the setting of metric spaces
\cite[Section 5]{PalmgrenLocalicCompletion}.

The following lemma corresponds to \cite[Theorem
5.2]{PalmgrenLocalicCompletion}. Note that the assumption in
\cite[Theorem 5.2]{PalmgrenLocalicCompletion} that $X$ is locally
compact is not necessary. 
\begin{lem}
  If $f \colon (X,M) \to (Y,M)$ is uniformly continuous on each open ball of $X$,
  then $r_{f}$ is a formal topology map from $\mathcal{U}(X)$ to
  $\mathcal{U}(Y)$.
\end{lem}
\proof
  We check \ref{FTM1}, \ref{FTM2}, \ref{FTM3a}, and \ref{FTM3b}.

 \ref{FTM1} Let $a \in U_X$. Choose any $d \in N$ and $\varepsilon
 \in \PRat$. Since $f$ is uniformly continuous on $a_{*}$, there
 exist $\rho \in M$ and $\delta \in \PRat$ such that
  \begin{align*}
    \left(\forall x,x' \in a_{*} \right) \rho(x,x') < \delta
    \implies d(f(x), f(x'))< \varepsilon.
  \end{align*}
  Then by  \ref{U2}, we have $a \cov_{X} a \downarrow
  \mathcal{C}^{\delta}_{\rho} \subseteq
  {r_{f}}^{-}\mathcal{C}^{2\varepsilon}_{d} \subseteq {r_{f}}^{-} U_Y$.

 \ref{FTM2} Let $b,c \in U_Y$ and $a \in {r_{f}}^{-}b \downarrow
 {r_{f}}^{-}c$. Then, there exist $b' <_Y b$ and $c' <_Y c$ such that
 $f[a_{*}] \subseteq b'_{*} \cap c'_{*}$. Write $b' =
 \ball_{d_{1}}(y,\delta)$ and $c' = \ball_{d_{2}}(z,\xi)$,
 and put $d = \sup\left\{ d_{1}, d_{2} \right\}$.
 Choose $\theta \in \PRat$ such that 
 $\ball_{d_1}(y,\delta + 2\theta) <_Y b$ and 
 $\ball_{d_2}(z,\xi + 2\theta) <_Y c$.
 Since $f$ is uniformly continuous on $a_{*}$, there exist
 $\rho \in M$ and $\varepsilon \in \PRat$ such that
  \begin{align*}
    \left(\forall x,x' \in a_{*} \right) \rho(x,x') < \varepsilon
    \implies d(f(x), f(x'))< \theta.
  \end{align*}
 Let $\ball_{\rho'}(x',\varepsilon') \in
 a \downarrow \mathcal{C}^{\varepsilon}_{\rho}$. Then,
 $f[\ball_{\rho'}(x',\varepsilon')_{*}] \subseteq \ball_{d}(f(x'),
 \theta)_{*}$. Since $\ball_{d}(f(x'), 2\theta) \in
 \ball_{d_{1}}(y,\delta + 2\theta) \downarrow \ball_{d_{2}}(z,\xi +
 2\theta) \subseteq b \downarrow c$,
 we have $\ball_{\rho'}(x',\varepsilon') \in {r_{f}}^{-}( b \downarrow c)$.
 Hence by \ref{U2}, we have $a \cov_X {r_{f}}^{-}(b \downarrow c)$.

 \ref{FTM3a} Obvious.

 \ref{FTM3b} For \ref{U1}, we have ${r_{f}}^{-} b \cov_{X}
 {r_{f}}^{-}\left\{ b' \in U_Y \mid b' <_{Y} b \right\}$  for all
 $ b\in U_Y$ by Lemma \ref{lem:misc} \eqref{lem:misc2}.
 For \ref{U2}, the argument is similar to the proof of the case \ref{FTM1}.
\qed

\begin{lem}\label{lem:USpaPtIsFunctor}
  Let $X$ be a locally compact uniform space, and let $Y$ be a
  complete uniform space. For any formal topology map $r \colon
  \mathcal{U}(X) \to \mathcal{U}(Y)$, the composition 
  \[
  f = {i_Y}^{-1} \circ \Pt(r) \circ i_{X}
  \]
  is uniformly continuous on each open ball of $X$. 
\end{lem}
\proof
See \cite[Theorem 5.1]{PalmgrenLocalicCompletion}.
\qed

\begin{lem}
  \label{lem:Faithful}
  Let $X$ and $Y$ be  complete uniform spaces, and let $f \colon X
  \to Y$ be a function which is uniformly continuous on each open ball of $X$.
  Then, the following diagram commutes.
  \[
  \xymatrix{
  X \ar[r]^-{i_{X}} \ar[d]_{f} & \Pt(\mathcal{U}(X))
  \ar[d]^{\Pt(r_{f})} \\
  Y   & \Pt(\mathcal{U}(Y)) \ar[l]_{{i_{Y}}^{-1}}
  }
  \]
\end{lem}
\proof
See \cite[Lemma 5.7 (i)]{PalmgrenLocalicCompletion}.
\qed

\begin{lem}
  \label{lem:Full}
  Let $X$ be a locally compact uniform space, and let $Y$ be a complete
  uniform space. Then, for any formal topology map $r \colon
  \mathcal{U}(X) \to \mathcal{U}(Y)$, we have 
  $r_{f} = r$,
  where $f \defeql {i_{Y}}^{-1} \circ \Pt(r) \circ i_{X}$.
\end{lem}
\proof
See \cite[Lemma 5.7 (ii)]{PalmgrenLocalicCompletion}.
\qed

\begin{lem}
  \label{lem:Functorial}
  Let $f \colon X \to Y$ and $g \colon Y \to Z$ be continuous functions between
  locally compact uniform spaces.  Then,  
  $
  r_{g\circ f} = r_{g} \circ r_{f}.
  $
\end{lem}
\proof
See the proof of \cite[Theorem 5.8]{PalmgrenLocalicCompletion}.
\qed

Similarly, we can show that $r_{\id_{X}} = \id_{\mathcal{U}(X)}$ for any
locally compact uniform space $X$. Hence, we conclude as follows.
\begin{thm}
  \label{thm:LocalicEmbedding}
  The localic completion $\mathcal{U}$ extends to a full and faithful
  functor
  \begin{align*}
    \mathcal{U} \colon \LKUSpa \to \OLKReg
  \end{align*}
  from the category
  of locally  compact uniform
  spaces $\LKUSpa$ to that of overt locally compact completely regular formal topologies
  $\OLKReg$.
\end{thm}
\proof
  For each morphism $f \colon X \to Y$ of $\LKUSpa$, define
  $\mathcal{U}(f) \defeql r_{f}$. Then, by Lemma \ref{lem:Functorial},
  $\mathcal{U}$ is a functor.  By Lemma \ref{lem:Faithful},
  $\mathcal{U}$ is faithful, and by Lemma \ref{lem:USpaPtIsFunctor}
  and Lemma \ref{lem:Full}, $\mathcal{U}$ is full.
\qed

By an abuse of terminology, we call the functor $\mathcal{U} \colon \LKUSpa
\to \OLKReg$ the localic completion of locally compact uniform spaces.

\subsection{Preservation of products}
Palmgren \cite{PalmgrenLocalicCompletion} showed that the localic
completion of locally compact metric spaces preserve finite products. We extend his result to the setting of
uniform spaces.
\begin{lem}
A binary product of locally compact uniform spaces (as generalised
uniform spaces) is locally compact.
\end{lem}
\proof
The proof is analogous to the metric case
\cite[Chapter 4, Proposition 12]{Bishop-67}.
\qed

\begin{lem}
  Let $(X,M)$ and $(Y,N)$ be generalised uniform spaces, and let
  $f \colon X \to Y$ be a homomorphism. Let $r$ and $s$ be the
  relations between $\mathcal{U}(X)$ and $\mathcal{U}(Y)$ as given by
  \eqref{def:UMap} and \eqref{eq:GUSFTop} respectively. Then, $r$
  and $s$ are equal as formal topology maps. 
\end{lem}
\proof
We must show that 
    $
    \sat r^{-} \ball_{\rho}(y,\delta)   = \sat s^{-} \ball_{\rho}(y,\delta) 
    $
for each $\ball_{\rho}(y,\delta) \in U_{Y}$.

First, suppose that  $\ball_{d}(x,\varepsilon) \mathrel{r}
\ball_{\rho}(y,\delta)$.  Then, there exists
$\ball_{\rho'}(y',\delta') <_{Y} \ball_{\rho}(y,\delta)$ such that
  $f[\Ball_{d}(x,\varepsilon)] \subseteq \Ball_{\rho'}(y',\delta')$.
  Choose $\theta \in \PRat$ such that $\rho(y,y') + \delta' + \theta <
  \delta$, and let $d' \in M$ such that $d'
  \mathrel{\omega_{f}} \rho'$. By
  $\ref{U2}$, we have
  \[
    \ball_{d}(x,\varepsilon) \cov_{X} \ball_{d}(x,\varepsilon) \downarrow
    \mathcal{C}^{\theta}_{d'}.
  \] 
  Let $\ball_{d^{*}}(x',\varepsilon') \in  \ball_{d}(x,\varepsilon)
  \downarrow \mathcal{C}^{\theta}_{d'}$. Then, $d^{*}
  \mathrel{\omega_{f}} \rho$, and
  we have 
  \[
    \rho(y,f(x')) + \varepsilon' \leq 
\rho(y,y') + \rho'(y', f(x')) + \varepsilon' < \rho(y,y') + \delta'
+ \varepsilon' < \delta.
  \]
  Thus, $\ball_{\rho}(f(x'),\varepsilon') <_{Y} \ball_{\rho}(y,\delta)$,
  so $\ball_{d}(x,\varepsilon) \cov_{X} s^{-} 
    \ball_{\rho}(y,\delta)$.

Next, suppose that
$\ball_{d}(x,\varepsilon) \mathrel{s} \ball_{\rho}(y,\delta)$.
Then, 
\[
  \rho(f(x),f(x')) \leq d(x,x')
\]
for any $x' \in X$, so $f[\Ball_{d}(x,\varepsilon)] \subseteq
\Ball_{\rho}(f(x),\varepsilon)$. Thus, 
$\ball_{d}(x,\varepsilon) \mathrel{r}\ball_{\rho}(y,\delta)$.
\qed

Since the projections from $X \times Y$ to $X$ and $Y$ are
homomorphisms, we have the following by Proposition \ref{prop:Prod}.
\begin{thm}\label{thm:Prod}
  For any locally compact uniform spaces $(X,M)$ and $(Y,N)$, we have
  \begin{align*}
    \mathcal{U}(X) \times \mathcal{U}(Y) \cong \mathcal{U}(X \times Y).
  \end{align*}
\end{thm}

Similarly, by Theorem \ref{prop:OmegaProd}, the localic completion
preserves inhabited countable products of compact uniform spaces. Note
that inhabited countable products of compact uniform spaces are again
compact since countable products of complete uniform spaces are
complete and inhabited countable products of total bounded uniform
spaces are total bounded (see \cite[Chapter 4, Problems
26]{Bishop-Bridges-85}). 

\begin{example}\label{eg:Cube}
  The Hilbert cube $\prod_{n \in \Nat} [0,1]$ is an inhabited
  countable product of the unit interval $[0,1]$, which is a compact
  metric space. The localic completion of $[0,1]$ is the formal unit
  interval $\UInt$, which is the overt weakly closed subtopology of the
  formal real $\FR$ (cf.\ Example \ref{eg:FR}) determined by the splitting
  subset
  \[
    \mathrm{Pos}_{\UInt} \defeql \left\{ (p,q) \in S_{\FR} \mid p < 1
  \amp 0 < q \right\}.
  \]
  Note that $\UInt$ is obtained as the localic completion of the
  rational interval $[0,1] \cap \Rat$ since the unit interval arises
  as its metric completion (cf.\ Proposition \ref{lem:DenseIso}).  In
  this case, we have $\mathcal{U}(\prod_{n \in \Nat} [0,1]) \cong
  \prod_{n \in \Nat}\UInt$.
\end{example}

\section{Connection to point-free completion}\label{sec:PFComp}
In [28],
Fox introduced the notion of a uniform formal topology, a formal
topology equipped with a covering uniformly, and established an
adjunction between the category of uniform spaces equipped with
covering uniformities and that of uniform formal topologies. He also
defined the completion of a uniform formal topology, which is classically equivalent to the
completion of a uniform locale by K\v{r}\'{i}\v{z}
\cite{KrizUniformCompletion}.

The construction given in this section on symmetric gus's is equivalent to applying the left adjoint of the
adjunction defined by Fox followed by the completion of uniform formal
topologies. We show that this construction is equivalent to the
localic completion.

First, we recall the relevant  notions from \cite[Chapter 6,
Section 2]{Fox05}. Our presentation is slightly different, but
equivalent to the one given in \cite{Fox05}.

\begin{defi}
Let $\mathcal{S}$ be an overt formal topology with positivity
$\Pos$. A \emph{cover} of
$\mathcal{S}$ is a subset $C \subseteq S$ such that $S \cov C$.
For covers $C,C' \in \Pow(S)$ of $\mathcal{S}$, define
\begin{align*}
  C \leqslant C' &\defeqiv \left( \forall a \in C \right) \left( 
  \exists a' \in C'\right) a \cov a', \\
  C <^{*} C' &\defeqiv \left(  \forall a \in C\right) \left(
  \exists a' \in C'\right) St_{C}(a) \cov a',
\end{align*}
where
$
St_{C}(a) \defeql \left\{ b \in C \mid \Pos(b \downarrow a)
\right\}.
$
A \emph{uniformity} on an overt formal topology $\mathcal{S}$
is a set $\mathcal{C}$ of covers of $\mathcal{S}$ such that
\begin{enumerate}
  \item  $\left( \forall C_1, C_2 \in \mathcal{C} \right)\left(
    \exists C_{3} \in \mathcal{C} \right) C_3 \leqslant C_1 \amp C_3 \leqslant  C_2$,
  \item  $\left( \forall C \in \mathcal{C} \right)\left( \exists C'
    \in \mathcal{C} \right) C' <^{*} C$,
  \item  $\left( \forall a \in S \right) a
    \cov \uc(a)$,
\end{enumerate}
where 
$
\uc(a) \defeql \left\{ b \in S \mid \left( \exists C \in
\mathcal{C} \right) St_{C}(b) \cov a 
\right\}.
$

A \emph{uniform formal topology} is a pair $(\mathcal{S},
\mathcal{C})$ where $\mathcal{S}$ is an overt formal topology
and $\mathcal{C}$ is a uniformity on $\mathcal{S}$.
\end{defi}

\begin{defi}
Let $(\mathcal{S},{\mathcal{C}})$ be a uniform formal topology.
Define a preorder $\leq$ on $S$ by $a \leq b \defeqiv a \cov b$.
The \emph{completion} of $\mathcal{S}$ is a formal topology
$\overline{\mathcal{S}} = (S, \bar{\cov}, \leq)$ 
inductively generated by the following axioms.
\begin{enumerate}[label=({CP}\arabic*)]
  \item  $a \bar{\cov}\uc(a)$;
  \item  $a \bar{\cov} C$ for each $C \in \mathcal{C}$;
  \item  $a \bar{\cov} \left\{ a \mid \Pos(a) \right\}$
\end{enumerate}
where 
$\Pos$ is the positivity of $\mathcal{S}$.
\end{defi}

A symmetric gus $(X,M)$ determines a uniform formal topology
$(\mathcal{S}_X,\mathcal{C}_M)$: the formal topology
$\mathcal{S}_X = (U_X,\cov_{X}^{*},\preceq_{X})$ is the usual topology induced
by the uniformity $M$, i.e.\
\begin{align*}
  U_X
  &\defeql \Ent(M) \times X,\\
  \ball_{d}(x,\varepsilon) \preceq_{X} \ball_{\rho}(y,\delta)
  &\defeqiv \ball_{d}(x,\varepsilon)_{*} \subseteq
  \ball_{\rho}(y,\delta)_{*}
  \iff \Ball_{d}(x,\varepsilon) \subseteq \Ball_{\rho}(y,\delta),\\
  \ball_{d}(x,\varepsilon) \cov_{X}^{*} U
  &\defeqiv \ball_{d}(x,\varepsilon)_{*} \subseteq
  U_{*},
\end{align*}
where we use the same notation for the elements of $U_X$ adopted in
Section \ref{subsec:LocComp}. Note that the positivity of $\mathcal{S}_{X}$ is $U_X$. The
uniformity $\mathcal{C}_{M}$ is given by 
\[
  \mathcal{C}_{M}
  \defeql  \left\{ \mathcal{C}^{\varepsilon}_{d}
  \mid (d, \varepsilon) \in \Ent(M) \right\},
\]
where $\mathcal{C}^{\varepsilon}_{d}$ is defined by \eqref{eq:UniCov}.
The pair $(\mathcal{S}_X,\mathcal{C}_M)$ is the
standard uniform formal topology with the usual topology
induced by the uniformity $M$.
The completion of $(\mathcal{S}_X,\mathcal{C}_M)$
is an inductively generated formal topology $\overline{\mathcal{S}_X}
= (U_X, \bcov_{X},\preceq_X)$, where $(U_X, \preceq_X)$ is the
underlying preorder of $\mathcal{S}_{X}$ and the cover $\bcov_{X}$ is
inductively generated by the following axioms:
\begin{enumerate}[label=({C}\arabic*)]
  \item\label{C1} $a \bcov_{X} \left\{ b \in U_X \mid b \prec_X a \right\}$;
  \item\label{C2} $a \bcov_{X} \mathcal{C}_{d}^{\varepsilon}\,$ for
    each $(d, \varepsilon) \in \Ent(M)$,
\end{enumerate}
for each $a \in U_X$, where
\[
  a \prec_{X} b \defeqiv \left( \exists (d,\varepsilon) \in \Ent(M) \right)
  \left( \forall c \in  \mathcal{C}^{\varepsilon}_{d} \right)
  c_{*} \meets a_{*} \imp c \preceq_{X} b \label{form:CovStrOrder}.
\]
The axioms \ref{C2} are equivalent to the following
axioms:
\begin{enumerate}[label=({C}\arabic*'),start=2]
  \item\label{C2'} $a \bcov_{X} \mathcal{C}_{d}^{\varepsilon}
    \downarrow_{\preceq_X} a \,$ for each $(d,\varepsilon) \in
    \Ent(M)$,
\end{enumerate}
where $\downarrow_{\preceq_X}$ is defined with respect to $\preceq_X$. 
Thus, the axioms \ref{C1}  and \ref{C2'} form a localised axiom-set on
$U_X$ with respect to  $\preceq_X$.

For each $(d,\varepsilon) \in \Ent(M)$, define a relation
$\prec_X^{(d,\varepsilon)}$ on $U_X$ by
\[
  a \prec_{X}^{(d,\varepsilon)} b \defeqiv 
  \left( \forall c \in  \mathcal{C}^{\varepsilon}_{d} \right)
  c_{*} \meets a_{*} \imp c \preceq_{X} b.
\]
Note that $a \prec_X b \iff \left( \exists (d,\varepsilon) \in
\Ent(M) \right) a \prec_X^{(d,\varepsilon)} b$.

\begin{lem}\label{lem:CXPos}
  Let $(X,M)$ be a symmetric gus. Then,
\begin{enumerate}
  \item $a \leq_X b \implies a
         \preceq_{X} b$,
  \item $a <_X b \implies a \prec_{X} b$,
\end{enumerate}
  for all $a, b \in U_X$.
\end{lem}
\proof
  (1) This is equivalent to Lemma \ref{lem:misc} \eqref{lem:misc3}.

  (2)  Let $a,b \in U_X$, and suppose that $a <_X b$.
   Write $a = \ball_{d}(x,\varepsilon)$ and $ b = \ball_{\rho}(y,\delta)$,
   and choose $\theta \in \PRat$ such that 
     $
     \rho(y,x) + \varepsilon  +  2 \theta < \delta.
     $
   We show that
    $
    a \prec_{X}^{(d,\theta)} b.
    $
    Let $c = \ball_{d}(z,\theta) \in \mathcal{C}^{\theta}_{d}$,
    and suppose that $c_{*} \meets a_{*}$. 
    Then, $d(x,z) < \varepsilon + \theta$, so 
    \begin{align*}
      \rho(y,z) + \theta
      &\leq \rho(y,x) + \rho(x,z) + \theta \\
      &\leq \rho(y,x) + d(x,z) + \theta \\
      &< \rho(y,x) + \varepsilon + 2\theta < \delta.
    \end{align*}
    Thus, $c \leq_{X} b$, so that $c \preceq_{X} b$. Hence
    $a \prec_{X}^{(d,\theta)} b$, and therefore $a \prec_{X} b$.
\qed

The following is a corollary of Lemma \ref{lem:Sqeq}, which holds for
any symmetric gus as well.
\begin{lem}\label{lem:UVCondIso}
  For any symmetric gus $(X,M)$, we have
  \begin{align*}
    a_{*} \subseteq b_{*} \amp b <_X c \implies a \cov_X c
  \end{align*}
  for all $a,b,c \in U_X$.
\end{lem}

\begin{thm}\label{thm:EquivLCompPFUComp}
  For any symmetric gus $(X,M)$, we have
    $
    \overline{\mathcal{S}_{X}} \cong \mathcal{U}(X).
    $
  That is, the localic completion of a symmetric gus is the point-free completion
  of the standard uniform formal topology induced by $M$.
  In particular, the localic completion of a symmetric gus is complete
  as a uniform formal topology.
\end{thm}
\proof
We define a binary relation $r_{X}$ on $U_X$ by
\[
  a \mathrel{r_X} b \defeqiv \left( \exists b' <_X b \right) a_{*} \subseteq
  b_{*}',
\]
and show that it is a surjective embedding from
$\overline{\mathcal{S}}_{X}$ to $\mathcal{U}(X)$.
\begin{enumerate}
  \item $r_{X}$ is a formal topology map:
  We check \ref{FTM1}, \ref{FTM2}, \ref{FTM3a}, and \ref{FTM3b}.

   \ref{FTM1} Let $a = \ball_d(x,\varepsilon) \in
  U_X$.  Then, we have $a \mathrel{r_{X}}
  \ball_{d}(x, 2\varepsilon)$, from which \ref{FTM1} follows.

  \ref{FTM2} Let $b,c \in U_X$, and let $a \in r_{X}^{-} b \downarrow
  r_{X}^{-}c$. Then, there exist $b' <_X b$ and $c' <_X c$ such that
  $a_{*} \subseteq b_{*}' \cap c_{*}'$.
  Write $b' = \ball_{d_1}(x, \varepsilon)$ and $c' =
  \ball_{d_2}(y,\delta)$, and choose $\theta \in \PRat$ such that
    $
    \ball_{d_1}(x,\varepsilon + 3\theta) <_X b$, and 
    $
    \ball_{d_2}(y,\delta+ 3\theta) <_X c.
    $
  Put $\rho = \sup\left\{ d_1,d_2 \right\}$, and
  let $a' = \ball_{\rho'}(z,\xi) \in \mathcal{C}^{\theta}_{\rho}
  \downarrow a$.
  Then, $a'_{*} \subseteq \ball_{\rho}(z,2\theta)_{*}$.
  Since $z \in a_{*}$, we have $d(x,z) + 3 \theta <  \varepsilon + 3\theta$, so that
  $\ball_{\rho}(z,3\theta) <_X \ball_{d_1}(x,\varepsilon + 3\theta)$.
  Similarly, we have
  $\ball_{\rho}(z,3\theta) <_X \ball_{d_2}(y,\delta + 3\theta)$.
  Hence, $a' \mathrel{r_{X}} \ball_{\rho}(z,3\theta)$ and
  $\ball_{\rho}(z,3\theta) \in b \downarrow c$. Therefore,
  $a \bcov_{X} r_{X}^{-}(b \downarrow c)$ by \ref{C2}.

  \ref{FTM3a} By Lemma \ref{lem:misc} \eqref{lem:misc1}.

  \ref{FTM3b} Preservation of the axiom \ref{U1} follows from Lemma
  \ref{lem:misc} \eqref{lem:misc2}. For \ref{U2}, let $(d,\varepsilon) \in
  \Ent(M)$. 
  Putting $\delta = \varepsilon/2$, we have
  $U_{X} \bcov_X
  \mathcal{C}^{\delta}_{d}
  \subseteq r_{X}^{-}\mathcal{C}^{\varepsilon}_{d}$ by \ref{C2}.

  \item $r_{X}$ is an embedding: We must show that 
    $
    a \bcov_{X} r_{X}^{-}r_{X}^{-*}\mathcal{A}_{X} \left\{ a
    \right\}
    $
  for all $a \in U_X$. Let $a \in U_X$,
  and let $a' \in U_X$ such that $a' \prec_{X} a$.
  Then, there exists $(d,\varepsilon) \in \Ent(M)$ such that 
  $a' \prec^{(d,\varepsilon)}_{X} a$.
  Choose $\theta \in \PRat$ such that $\theta < \varepsilon$,
  and let $b \in \mathcal{C}^{\theta}_{d} \downarrow a'$.
  Then, there exists $\ball_{d}(x,\theta) \in
  \mathcal{C}^{\theta}_{d}$
  such that $b_{*} \subseteq \ball_{d}(x,\theta)_{*}$, and thus
  $b \mathrel{r_{X}} \ball_{d}(x,\varepsilon)$.
  Let $b' \in r_{X}^{-}\ball_{d}(x,\varepsilon)$.
  Then, we have $b'_{*} \subseteq \ball_{d}(x,\varepsilon)_{*}$.
  Since $a'_{*} \meets \ball_{d}(x,\varepsilon)_{*}$, 
  we have 
  $\ball_{d}(x,\varepsilon)_{*} \subseteq a_{*}$.
  Hence  $b' \preceq_{X} a$, so
  $\ball_{d}(x,\varepsilon) \in r_{X}^{-*} \mathcal{A}_{X}
  \left\{ a \right\}$.
  Therefore, by \ref{C1} and  \ref{C2}, we have
  $a \bcov_{X} r_{X}^{-} r_{X}^{-*} \mathcal{A}_{X} \left\{ a \right\}$,
  as required.

  \item $r_{X}$ is a surjection:
  We must show that
  \[
  r_{X}^{-}a \bcov_{X} r_{X}^{-}U \implies a \cov_X U
  \]
  for all $a\in U_X$ and $U \subseteq U_X$.  Since
    $
    b <_X a \implies b \mathrel{r_{X}} a
    $
    for all $a, b \in U_X$, it suffices to show that
    \[
    a \bcov_{X}
    r_{X}^{-}U \implies a \cov_X U
    \]
  for all $a \in U_X$ and $U \subseteq U_X$ by \ref{U1}.
  Given $U \subseteq U_X$, define a predicate $\Phi_U$ on
  $U_X$ by
  \[
    \Phi_U(a) \defeqiv \left( \forall b \in U_X \right) b_{*} \subseteq
    a_{*} \implies b \cov_{X} U.
  \]
  Then, it suffices to show that
  \[
  a \bcov_{X} r_{X}^{-}U \implies \Phi_U(a)
  \]
  for all $a \in U_X$.
  This is proved by induction on $\bcov_{X}$. We must
  check the conditions \ref{ID1} -- \ref{ID3}.

  \ref{ID1} Suppose that $a \in r_{X}^{-}U$, and let $b \in U_X$ such that
  $b_{*} \subseteq a_{*}$. Then, there exist $c \in U$ and $c' <_X c$
  such that $a_{*} \subseteq c'_{*}$. Thus $b_{*} \subseteq
  c'_{*}$, and hence, $b \cov_{X} c$ by Lemma
  \ref{lem:UVCondIso}. Therefore $\Phi_U(a)$.

  \ref{ID2} This directly follows from the definition of $\preceq_X$.

  \ref{ID3} We need to check the axioms \ref{C1} and \ref{C2'}.

  \ref{C1} $\displaystyle{\frac{\left( \forall b  \prec_{X} a 
   \right)\Phi_U(b)}{\Phi_U(a)}}$:
  Suppose that $\Phi_U(b)$ holds for all $b \prec_{X}
  a$.
  Let $c \in U_X$ such that $c_{*} \subseteq a_{*}$, and let  $b \in
  U_X$ such that $b <_X c$. 
  Then, $b \prec_{X} c$ 
  by Lemma \ref{lem:CXPos}, and so $b
  \prec_{X} a$.
  Thus, $\Phi_U(b)$, and so $b \cov_{X} U$.
  Hence $c \cov_{X} U$ by \ref{U1}. Therefore $\Phi_U(a)$.

  \ref{C2'} $\displaystyle{\frac{\left( \forall b \in
    \mathcal{C}^{\varepsilon}_{d} \downarrow_{\preceq_X} a \right)\Phi_U(b)}{\Phi_U(a)}}$ for
  each $(d,\varepsilon) \in \Ent(M)$:
  Let $(d,\varepsilon) \in \Ent(M)$, and 
  suppose that $\Phi_U(b)$ holds for all $b \in \mathcal{C}^{\varepsilon}_{d} \downarrow_{\preceq_X} a$.
  Let $c \in U_X$ such that  $c_{*} \subseteq a_{*}$.
  Let $b \in \mathcal{C}^{\varepsilon}_{d} \downarrow c$.
  Then, $b \in \mathcal{C}^{\varepsilon}_{d} \downarrow_{\preceq_X} a$,
  and so $\Phi_U(b)$. Thus, $b \cov_{X} U$, and hence $c \cov_{X} U$ by
  \ref{U2}. Therefore $\Phi_U(a)$. \qedhere
\end{enumerate}

\section{Further work}
  Curi \cite{CuriCollectionOfPoints} developed a theory of uniform
  formal topologies that is different from that of Fox \cite{Fox05}.  In
  the same paper, Curi sketched another embedding from the category
  of uniform spaces\footnote{Curi's notion of uniform space is that of
  symmetric generalised uniform space in our terminology.} and
  uniformly continuous functions into uniform formal topologies.

  Curi's embedding, however, sends a uniform space to a formal topology
  which has a usual topology induced by the uniformity. Hence,
  without Fan theorem, his embedding cannot be shown to preserve
  compactness and local compactness of uniform spaces as our localic
  completion does. This is one of our motivations to develop the
  localic completion of generalised uniform spaces. However, a uniform
  formal topology that arises from Curi's embedding might be
  `completed' in a suitable sense to give a topology
  equivalent to the localic completion. Development of completions of
  uniform formal topologies in Curi's sense and its comparison to our
  localic completion and Fox's notion of uniform formal topology are
  left for the future.

\section*{Acknowledgement}
The author is grateful to the anonymous referees for their 
generous suggestions which significantly improved the
presentation of the paper.
This work is supported by Core-to-Core Program A.~Advanced
Research Networks by Japan Society for the Promotion of Science
(JSPS).
\bibliographystyle{alpha}

\end{document}